\documentclass[a4paper,10pt]{article}

\usepackage{amsmath,amssymb,wasysym,amsfonts,amsthm,amsthm}

\usepackage{empheq,hyperref,url,enumerate}
\usepackage[a4paper,margin=1.0in]{geometry}

\newtheorem{result}{Result}

\theoremstyle{remark}
\newtheorem{rmk}[result]{Remark}

\numberwithin{equation}{section}
\numberwithin{result}{section}

\usepackage{color}
\usepackage[dvipsnames]{xcolor}

\usepackage{cite}

\usepackage{graphicx,epsfig,subfig}
\usepackage{tikz-cd}

\usepackage[british]{babel}
\usepackage{soul}

\newcommand{\eps}{\varepsilon}

\title{A PREDATOR--2 PREY FAST--SLOW DYNAMICAL SYSTEM FOR RAPID PREDATOR EVOLUTION} 

\author{Sofia H. Piltz \thanks{Department of Applied Mathematics and Computer Science and National Institute of Aquatic Resources, Technical University of Denmark, Asmussens all\'{e}, Bygning 303B, 2800 Kongens Lyngby, Denmark and Department of Mathematics, University of Michigan, 2074 East Hall, 530 Church Street, Ann Arbor, Michigan, 48109-1043, USA \mbox{(\href{mailto:piltz@umich.edu}{piltz@umich.edu})}.}
\and Frits Veerman \thanks{Wolfson Centre for Mathematical Biology, Mathematical Institute, University of Oxford, Andrew Wiles Building, Radcliffe Observatory Quarter, Woodstock Road, Oxford, OX2 6GG, UK \mbox{(\href{mailto:veerman@maths.ox.ac.uk}{veerman@maths.ox.ac.uk})}; School of Mathematics, University of Edinburgh, James Clerk Maxwell Building, Peter Guthrie Tate Road, Edinburgh, UK \mbox{(\href{mailto:f.veerman@ed.ac.uk}{f.veerman@ed.ac.uk})}.} 
\and Philip K. Maini \thanks{Wolfson Centre for Mathematical Biology, Mathematical Institute, University of Oxford, Andrew Wiles Building, Radcliffe Observatory Quarter, Woodstock Road, Oxford, OX2 6GG, UK and CABDyN Complexity Centre, University of Oxford, Oxford, OX1 1HP, UK \mbox{(\href{mailto:maini@maths.ox.ac.uk}{maini@maths.ox.ac.uk})}.}
\and Mason A. Porter\thanks{Oxford Centre for Industrial and Applied Mathematics, Mathematical Institute, University of Oxford, Andrew Wiles Building, Radcliffe Observatory Quarter, Woodstock Road, Oxford, OX2 6GG, UK; CABDyN Complexity Centre, University of Oxford, Oxford, OX1 1HP, UK; and Department of Mathematics, UCLA, Los Angeles, CA 90095, USA \mbox{(\href{mailto:mason@math.ucla.edu}{mason@math.ucla.edu})}.}}

\begin{document}
\maketitle

\begin{abstract}

We consider adaptive change of diet of a predator population that switches its feeding between two prey populations. We develop a novel 1 fast--3 slow dynamical system to describe the dynamics of the three populations amidst continuous but rapid evolution of the predator's diet choice. The two extremes at which the predator's diet is composed solely of one prey correspond to two branches of the three-branch critical manifold of the fast--slow system. By calculating the points at which there is a fast transition between these two feeding choices (i.e., branches of the critical manifold), we prove that the system has a two-parameter family of periodic orbits for sufficiently large separation of the time scales between the evolutionary and ecological dynamics. Using numerical simulations, we show that these periodic orbits exist, and that their phase difference and oscillation patterns persist, when ecological and evolutionary interactions occur on comparable time scales. Our model also exhibits periodic orbits that agree qualitatively with oscillation patterns observed in experimental studies of the coupling between rapid evolution and ecological interactions.

\end{abstract}




\section{Introduction}\label{section_introduction}

Organisms can adapt to changing environmental conditions---such as prey availability, predation risk, or temperature---by changing their behavior. For example, in \emph{prey switching}, a predator changes its diet or habitat in response to prey abundances. This is an example of \emph{phenotypic plasticity} \cite{Kellyetal_phenplast2012}, in which the same genotype can express different phenotypes in different environments. However, adaptivity can also be expressed as an evolutionary change in \emph{traits} (i.e., properties that affect how well an individual performs as an organism \cite{McGilletal2006}) via genomic changes of a predator and/or prey \cite{Fussman2007}. If such evolution occurs on a time scale of about 1000 generations and can be observed in laboratory conditions, it is construed to be a `rapid' evolutionary change of traits \cite{Fussman2007}. Rapid evolutionary changes have been observed in a wide variety of organisms, ranging from mammals \cite{Pelletier2007} to bacteria \cite{BohannanLenski2000}, and both in predators (e.g., in traits that involve resource consumption \cite{GrantGrant2002} or the ability to counteract prey defense mechanisms \cite{Hairstonetal1999}) and in prey (e.g., in traits that involve predator avoidance \cite{Jones2009,Yoshidaetal2003}). 
Understanding the dependencies between rapid evolution and ecological interactions is fundamental for making accurate predictions of a population's ability to adapt to, and persist under, changing environmental conditions \cite{Schoener2011,Ellneretal2011}. For example, rapid evolutionary change of traits has been observed in a plankton predator--prey system \cite{Fussmann2000,Yoshidaetal2003}, which is a good example system for studying the coupling between rapid evolution and predator--prey interaction due to its short generation times and the tractability of genetic studies of it \cite{Jones2009}.

Adaptive change of feeding behavior can be incorporated into a dynamical system of predator--prey interaction in multiple ways. For example, one can represent prey switching with a Holling type-III functional response \cite{Holling1965}, consider the densities of different prey as system variables \cite{AbramsMatsuda2003,Postetal2000}, or use information on which prey type was last consumed \cite{vanLeeuwen2007,vanLeeuwen2013}. Such formulations lead to smooth dynamical systems, but one can also model a predator that switches prey using a \emph{piecewise-smooth dynamical system} \cite{pw-sBook,pws-scholarpedia} in which continuous temporal evolution of predator and prey populations alternates with abrupt events that correspond to points at which the predator changes its diet or habitat \cite{Krivan1996,Krivan1997}. 


\paragraph{Rapid evolution.}

Several theoretical and empirical investigations have considered the effect of rapid evolutionary change of traits on predator--prey dynamics (see \cite{Fussman2007} for a review), including examples in which ecological and evolutionary dynamics have been assumed to occur on (1) comparable time scales or (2) disparate time scales. An example of (1) is the occurrence of out-of-phase cycles between small zooplankton (i.e., a predator) and genetically variable clonal lines of algae (i.e., prey) populations observed in the experiments in \cite{Fussmann2000,Yoshidaetal2003}, which were reproduced using a mathematical model with contemporaneous evolutionary and ecological dynamics \cite{JonesEllner2007}. The model in \cite{JonesEllner2007} suggests that the cycles emerge from prey evolution, especially when there is a small (energy) cost associated with the prey defense mechanism. Examples of (2) include situations in which evolutionary change occurs on either a slower \cite{KhibnikKondrashov1997} or on a faster \cite{CortezEllner2010} time scale than ecological interactions. Consequences that ecological dynamics can have on trait evolution have also been studied using the mathematical framework of \emph{adaptive dynamics} \cite{Geritzetal1998,ZuWang2013}, where evolution is assumed to occur on a slower time scale than ecological interactions. In the present paper, we aim to provide insight on how the evolution of traits arises in population dynamics, and we thus concentrate on studying the limit in which trait evolution occurs on a much faster time scale than predator--prey interactions. When the time scales can be separated, one can use the framework of \emph{fast--slow dynamical systems} \cite{KuehnMTS2015} to introduce and exploit a time-scale separation between evolutionary and ecological dynamics to reduce the dimensionality of the system of equations that describe the evolutionary and ecological interactions. \textcolor{black}{For a short introduction to fast--slow dynamical systems, see Appendix \ref{appendix_gspt}.} 

When evolutionary change is faster than ecological interactions, the fast--slow dynamical system introduced in \cite{CortezEllner2010} can preserve the qualitative properties of dynamics in a predator--evolving-prey model in which ecological changes and evolution occur on the same time scale \cite{JonesEllner2007}. Additionally, similar to a model with only one time scale \cite{JonesEllner2007}, a fast--slow dynamical system with rapid prey evolution reproduces experimentally-observed out-of-phase predator-prey oscillations \cite{Yoshidaetal2003}. However, such oscillations are not present in the analogous model without rapid evolution \cite{CortezEllner2010}. By exploiting the time-scale separation between fast evolution and slow ecological changes, the general theory of either an evolving prey or an evolving predator \cite{CortezEllner2010} has been extended to cover the case in which both predator and prey evolve \cite{CortezTheorEcol2015}. There exist general conditions for determining which type of cyclic dynamics are possible in a system of coevolving prey and predator \cite{CortezTheorEcol2015}. Such dynamics involve cycles that exhibit (i) counterclockwise or clockwise orientation in the predator--prey phase plane, (ii) a half-phase difference between the predator and prey oscillations, and (iii) `cryptic' cycles in which the predator population cycles while the prey population is approximately constant. Interestingly, a situation in which both predator and prey evolve can generate clockwise cycles, which have been identified in empirical data sets from systems such as phage--cholera, mink--muskrat, and gyrfalcon--ptarmigan \cite{CortezWeitz2014}. This contrasts with traditional Lotka--Volterra predator-prey cycles, which have a counterclockwise orientation in the phase plane,
with a quarter-phase lag between the predator and prey oscillations \cite{Murray1993}.

\paragraph{Our approach.} 

In the present paper, we use a similar approach to \cite{CortezEllner2010,CortezTheorEcol2015} and develop a novel (to our knowledge) fast--slow dynamical system for a predator switching between two groups of prey species. In contrast to \cite{CortezEllner2010,CortezTheorEcol2015}, we make a simplifying assumption of unlimited prey growth (e.g., because of favorable environmental conditions) and use Lotka--Volterra functional responses between predator and prey. As a result, we can prove that there exists a family of periodic orbits in a system of one predator and two different prey species. 
As we discuss in Section \ref{section_discussion}, some of our orbits agree qualitatively with patterns observed in both laboratory experiments and field research. In addition to the potential utility of mathematical modeling (and using fast--slow systems) for understanding the coupling between ecology and evolution \cite{Fussman2007}, our motivation for constructing our model comes from our earlier work that suggests that adaptive prey switching of a predator is a possible mechanistic explanation for patterns observed in data on freshwater plankton \cite{OurPaper}. In the model in \cite{OurPaper}, the switch in the predator's feeding behavior is discontinuous. In the present paper, we relax this assumption and consider a rapid but continuous change in the predator's feeding choice.\footnote{In other work, we consider two types of regularizations of the discontinuous switch that do not introduce a time-scale difference into the model \cite{ourSmoothPaper}.}

There exists theory both for regularizing a given piecewise-smooth dynamical system to create a fast--slow dynamical system \cite{Sotomayoretal1996,Teixeira2012} and for approximating a fast--slow system using a piecewise-smooth system that preserves---both qualitatively and quantitatively---key characteristics (such as singularities and bifurcations) of the original fast--slow system \cite{DesrochesJeffreyFoldedNode2015}. However, we instead construct our model from a biological perspective using a concept from quantitative genetics. A clear understanding of the trajectories of solutions of a fast--slow system makes it possible to compare a fast--slow system as an `ecologically-obtained' regularization (which we construct using fitness-gradient dynamics \cite{Abramsetal1993,Lande1982}) of a piecewise-smooth system with a regularization obtained from a mathematical viewpoint (e.g., using a method based on a blow-up technique \cite{Sotomayoretal1996,Teixeira2012}). We aim to shed light on the differences and similarities that phenotypic plasticity and rapid evolution, as two different mechanisms, can generate in a model for the population dynamics of an adaptively feeding predator and its two prey. In Section \ref{section_discussion}, we compare and contrast the dynamics exhibited by the earlier piecewise-smooth model \cite{OurPaper} and the fast--slow system that we construct as an ecological regularization of it.

Putting aside our motivating applications in ecology, we note in passing that classical examples of fast--slow dynamical systems include the Van der Pol \cite{vanDerPolSystem1920} and Fitzhugh--Nagumo equations \cite{FitzHugh1961,Nagumoetal1962}. The former was used originally to describe the dynamics of an electrical circuit with an amplifying valve (and has subsequently been used for numerous other applications), and the latter is a simplified version of the Hodgkin--Huxley nerve-axon model \cite{HodgkinHuxley1952} from neuroscience \cite{termanbook}. There are several other applications of such multiple time-scale systems---including pattern formation \cite{GiererMeinhardt1972}, opening and closing of plant leaves \cite{Forterre2013}, ocean circulation \cite{KowalczykGlendinning2011}, critical transitions in climate change \cite{kuehn2011mathematical}, and more. We also note that multiple time-scale systems can be studied using several different techniques from singular perturbation theory \cite{Hinch1991,KuehnMTS2015}. Examples of such techniques include matched asymptotic expansions \cite{KevorkianCole1996,bo99} and geometric singular perturbation theory \cite{Jones1995,KuehnMTS2015}. Because we are interested in constructing periodic orbits and understanding their bifurcations in the fast--slow dynamical system for prey switching (see Section \ref{section_fastSlowModel}), we use the latter to analyze our model.


\paragraph{Outline of our paper.} 

The rest of our paper is organized as follows. In Section \ref{section_fastSlowModel}, we formulate a 1 fast--3 slow dynamical system for a predator feeding on two prey populations in the presence of rapid predator evolution. The three slow variables of the system correspond to the populations of the predator and the two prey. We model a predator trait that represents the predator's feeding choice as the fast variable of the system. This model construction allows us to use geometric singular perturbation theory to gain insight into the effects on population dynamics of an evolutionary change of a predator trait that occurs on a time scale that is comparable to that of the predator--prey interaction. In Section \ref{section_preliminaries}, we derive expressions for the critical manifold and the slow and fast subsystems. We then use these results in Section \ref{section_constructionOfThePeriodicOrbitSection} to explicitly construct periodic orbits that are exhibited by the 1 predator--2 prey fast--slow system. We obtain these expressions for the periodic orbits by studying the singular limit in which the ratio $\eps$ of the fast to the slow time scale goes to $0$. In Section \ref{section_QualAsp}, we highlight some ecologically relevant qualitative aspects of the constructed periodic orbits. We then use numerical continuation in Section \ref{section_numericalContinuation} to investigate how the periodic orbits persist for $\eps>0$ as we perturb the system. Finally, we discuss the findings of our study in Section \ref{section_discussion}. \textcolor{black}{We briefly review geometric singular perturbation theory in Appendix \ref{appendix_gspt}, and we give additional details about finding families of singular periodic orbits in Appendix \ref{appendix:findsols}. In supplementary information (SI), we provide Mathematica notebooks containing our numerical code for finding and visualizing periodic orbits. We also provide associated data files containing the results of our numerical computations.}


\section{The fast--slow 1 predator-2 prey model}
\label{section_fastSlowModel}

We begin our formulation of a 1 predator--2 prey fast--slow system in Section \ref{subsection_popDyn_fastslowmodel_derivation} by constructing an equation for the temporal evolution of a predator population ($z$) that adaptively changes its diet between two prey populations ($p_1$ and $p_2$). Our fast--slow model is based on four principal assumptions. We assume that the organisms (1) have a large population size, (2) live in a well-mixed environment, and (3) can be aggregated into groups of similar species. Consequently, we can represent the predator--prey interaction as a low-dimensional system of ordinary differential equations. We \textcolor{black}{presume that the predator--prey interaction is such that it is possible for evolutionary changes of traits to occur on a time scale that is comparable to that of ecological interaction}\textcolor{black}{. In previous work \cite{OurPaper}, this evolutionary change was modeled as an instantaneous switch. To bridge the gap between this model of instantaneous evolutionary change and the ecological presumption that ecological traits change on a time scale that is comparable to the ecological interaction of the species, we study the limit in which (4) a predator trait undergoes rapid evolution on a faster time scale than that of the population dynamics. This gives insight into contemporaneous demographic and evolutionary changes.} In Section \ref{section_evolDyn_fastslowmodel_derivation}, we define an expression for the temporal evolution of a predator trait ($q$) that represents the predator's desire to consume each prey.  \textcolor{black}{In Section \ref{section_numericalContinuation}, we examine possible insights into interactions between ecological and evolutionary dynamics when these occur on a comparable time scale.}


\subsection{Ecological dynamics}\label{subsection_popDyn_fastslowmodel_derivation}

We assume that the predator's desire $q$ to consume prey is bounded between its smallest and largest feasible values ($q_{S}$ and $q_{L}$, respectively). For simplicity, we consider exponential prey growth and a linear functional response between the predator growth and prey abundance {\cite{Murray1993}. We thereby obtain the following system of differential equations for the population dynamics of the 1 predator--2 prey fast--slow system:
\begin{align} \label{fastslow_popDynMoregeneral} 
	\frac{\text{d} p_1}{\text{d} t}&=r_1 p_1-(q-q_{S})\beta_1 p_1 z \,,\nonumber \\
	\frac{\text{d} p_2}{\text{d} t}&=r_2 p_2-(q_{L}-q)\beta_2 p_2 z  \,,\\ \nonumber
	\frac{\text{d} z}{\text{d} t}&=e(q-q_{S})\beta_1 p_1 z+e(q_{L}-q)q_2\beta_2 p_2 z - m z \,, 
\end{align}
where $r_1$ and $r_2$ (with $r_1, r_2>0$) are the respective per capita growth rates of prey $p_1$ and $p_2$, the parameters $\beta_1$ and $\beta_2$ are the respective death rates of the prey $p_1$ and $p_2$ due to predation, $e>0$ is the proportion of predation that goes into predator growth, $q_2 \in [0,1]$ is the nondimensional parameter that represents the extent of preference towards prey $p_2$, and $m>0$ is the predator's per capita death rate. One can also interpret the parameter $q_2$ as a factor that scales the benefit that the predator obtains from feeding on prey $p_2$.

For simplicity, we let $\beta_1=\beta_2$ (which we can take to be equal to $1$ by rescaling) to omit $\beta_1$ and $\beta_2$ in our calculations in Section \ref{section_preliminaries}. In doing so, we assume that the predator exhibits adaptive diet choice by adjusting its feeding choice (i.e., whether the predator is feeding on prey $p_1$ or on prey $p_2$) rather than its attack rate based on the prey densities. We also require that the extreme when $q$ is at its minimum (i.e., $q=q_{S}$) corresponds to the case in which the predator is feeding solely on prey $p_2$. Similarly, we require that the extreme when $q$ is at its maximum (i.e., $q=q_{L}$) corresponds to the case in which the predator feeds solely on prey $p_1$. We thereby assume that $q$ is bounded between $q_L$ and $q_S$. Without loss of generality, we choose $q_L = 0$ and $q_S = 1$. These assumptions simplify the system \eqref{fastslow_popDynMoregeneral}, which represents the fast--slow 1 predator-2 prey population dynamics, to
\begin{align} \label{fastslow_popDyn} 
	\frac{\text{d} p_1}{\text{d} t}&=r_1 p_1-q p_1 z \,,\nonumber \\
	\frac{\text{d} p_2}{\text{d} t}&=r_2 p_2-(1-q)p_2 z  \,,\\ \nonumber
	\frac{\text{d} z}{\text{d} t}&= e q p_1 z+e(1-q)q_2 p_2 z - m z \,. 
\end{align}


\subsection{Evolutionary dynamics}
\label{section_evolDyn_fastslowmodel_derivation}

We assume that the adaptive change in the predator's trait $q$ follows {\emph{fitness-gradient dynamics}. In other words, we assume that the rate of change of the mean trait value is proportional to the fitness gradient of an individual with this mean trait value \cite{Abramsetal1993}. Fitness-gradient dynamics was used for defining trait dynamics of rapid predator evolution in \cite{CortezEllner2010,CortezTheorEcol2015,CortezWeitz2014}, wherein fast--slow dynamical systems were proposed as a general framework for gaining insight into evolutionary and ecological dynamics that occur on a comparable time scale. In the original form of fitness-gradient dynamics in \cite{Abramsetal1993}, the fitness $F$ of an individual is assumed to be frequency-dependent. That is, $F=F(q^*,q)$, where $q^*$ is the trait value of an individual and $q$ is the mean trait value of the population. In the present paper, we determine fitness as the net per capita growth rate of the predator population. The rate of change of the mean population trait value is then governed by fitness-gradient dynamics as follows:
\begin{equation}
	\frac{\text{d} q}{\text{d} t}\propto\frac{\partial}{\partial q^*}\left.\left(\frac{1}{z}\frac{\text{d} z}{\text{d} t}(p_1,p_2,z,q^*,q)\right)\right\vert_{q^*=q}\,.\label{fitnessGradient_predator}
\end{equation}
Following ecological considerations, if adaptation occurs by genetic change, the rate constant that describes \eqref{fitnessGradient_predator} is the \emph{additive genetic variance} \cite{Lande1982} (i.e., genetic variance due to genes whose alleles contribute additively to the trait value). For simplicity, we assume that the fitness of an individual depends only on the mean trait value of the population and not on the distribution of the individual's trait. That is, $F(q,q^*)=F(q)$. \textcolor{black}{Note that the use of this simplification implies that we assume the distribution of the trait to be sufficiently narrow.} Additionally, we describe the additive genetic variance of the predator's desire to consume each prey type by a bounding function that limits the predator trait between its smallest ($q_{S}=0$) and largest ($q_{L}=1$) feasible values. Furthermore, by assuming that the predator trait evolves on a faster time scale than the population dynamics (where the separation of time scales is given by $\eps$), we see that the temporal evolution of the predator trait takes the following form:
\begin{align}
	\eps\frac{\text{d} q}{\text{d} t}&=(q-q_{S})(q_{L}-q)V\frac{\partial}{\partial q}\left(\frac{1}{z}\frac{\text{d} z}{\text{d} t}(p_1,p_2,z,q)\right) \nonumber \\
	&=q(1-q)Ve\left(p_1-q_2 p_2\right)\,, \label{q_timeEvolution_fastslowSystem}
\end{align}
where $V$ is a nondimensional constant and is part of the additive genetic variance term $q(1-q)V$.

When using fitness-gradient dynamics, we model evolutionary dynamics at the phenotypic level without incorporating detailed information about genotypic processes (e.g., principles of Mendelian inheritance \cite{Fussman2007}). Using this simplified approach, we can incorporate an equation for the predator trait directly into the system of 1 predator--2 prey dynamics and obtain an analytically tractable differential-equation model for coupled ecological and evolutionary dynamics. Because of its simplifying assumptions on genotypic processes and the laws of inheritance, fitness-gradient dynamics gives an incomplete understanding of interactions between ecological and evolutionary dynamics \cite{Fussman2007}. Nevertheless, there is evidence that fitness-gradient dynamics can still be an appropriate approximation for modeling evolutionary dynamics even when its simplifying assumptions do not hold \cite{Gomulkiewicz1998,Abramsetal1993}. 


\subsection{Coupled ecological and evolutionary dynamics}
\label{section_CoupledEcolEvolDynamics} 

By combining the ecological dynamics in \eqref{fastslow_popDyn} with the evolutionary dynamics in \eqref{q_timeEvolution_fastslowSystem}, we obtain the following fast--slow 1 predator--2 prey system with predator evolution: 
\begin{align}
	\frac{\text{d} p_1}{\text{d} t}&=\dot{p_1}=g_1(p_1,p_2,z,q)=r_1p_1-qp_1z \nonumber \,,\\
	\frac{\text{d} p_2}{\text{d} t}&=\dot{p_2}=g_2(p_1,p_2,z,q)=r_2p_2-(1-q)p_2z  \label{fastslowFullSystem} \,,\\ 
	\frac{\text{d} z}{\text{d} t}&=\dot{z}=g_3(p_1,p_2,z,q)=eqp_1z+e(1-q)q_2p_2z-mz \nonumber \,,\\
	\eps\frac{\text{d} q}{\text{d} t}&=\dot{q}=f(p_1,p_2,q)=q(1-q)Ve(p_1-q_2p_2)\,. \nonumber
\end{align} 
When $q=1$ in \eqref{fastslowFullSystem}, the predator feeds only on prey $p_1$. Likewise, when $q=0$, the vector field of the fast--slow system \eqref{fastslowFullSystem} corresponds to a situation in which the predator's diet is composed solely of prey $p_2$. Consequently, there is exponential growth in the population of the prey type that is not being 
preyed upon.

\section{Analytical setup}
\label{section_preliminaries}

\textcolor{black}{In this section, we use geometric singular perturbation theory to aid in the analysis of the 1 fast--3 slow model in \eqref{fastslowFullSystem}. See Appendix \ref{appendix_gspt} for a brief introduction to geometric singular perturbation theory.}


\subsection{Rescaling of the system \eqref{fastslowFullSystem}}\label{section_rescalingthesystem}

To keep our analysis as clear as possible, we rescale \textcolor{black}{the system \eqref{fastslowFullSystem}} to maximally reduce the number of parameters. Using the rescaling
\begin{equation}
	 t \to \frac{t}{r_1}, \;p_1 \to \frac{m\,r_1}{e} p_1, \;p_2 \to \frac{m\,r_1}{e\,q_2}p_2, \;z \to r_1\,z, \;m\to r_1\,m, \;r_2 \to r\, r_1, \;\eps \to \eps\, m\,V\,,
\end{equation}
we obtain
\begin{align}
  \dot{p_1} &= (1 - q\,z)\,p_1\,, \nonumber\\
  \dot{p_2} &= (r - (1-q)z)\,p_2\,, \nonumber\\
  \dot{z} &= (q\,p_1 + (1-q)p_2-1)\,m\,z\,, \label{rescaled_system}\\
  \eps \dot{q} &= q(1-q)\,(p_1 - p_2)\,,\nonumber
\end{align}
where $r$ and $m$ are free parameters. Without loss of generality, we can assume that $0<r<1$ (i.e., that $r_1 > r_2$).


\subsection{Linearization around the coexistence equilibrium of \eqref{rescaled_system}}\label{section_lin_centrecentre}

The system \eqref{rescaled_system} has a unique \textcolor{black}{coexistence steady state} at $(p_1,p_2,z,q) = \left(1,1,1+r,\frac{1}{1+r}\right)$. The local behavior near this equilibrium is characterized by the eigenvalues ($\lambda$) of the linearization of \eqref{rescaled_system} around it. These eigenvalues obey the characteristic equation
\begin{equation}\label{this}
	 \lambda^4 + \frac{m + 2 r + m r^2}{1+r}\,\lambda^2 + m r = 0\,.
\end{equation}
For all $m>0$ and $0<r<1$, equation \eqref{this} has four purely imaginary solutions, so the \textcolor{black}{coexistence} equilibrium is of center--center type. For an equilibrium of this type, local analysis alone is in general extremely complicated and intricate (see, e.g., Chapter 7.5 of \cite{GuckenheimerHolmes}). Moreover, any local periodic solution stays close to the value $q=\frac{1}{1+r}$, whereas our goal is to \textcolor {black}{investigate different types of periodic orbits that can arise from rapid predator evolution, which occurs when $q$ varies between---and comes close to---the extremal values $0$ and $1$.} We therefore abandon the standard linearization approach and turn to geometric singular perturbation theory to obtain far-from-equilibrium periodic orbits (see Section \ref{section_constructionOfThePeriodicOrbitSection}).


\subsection{Analysis of the system \eqref{rescaled_system} as a fast--slow system}
\label{section_analysis_of_theFastSlowSystem}

In this section, we specify the fast and slow subsystems of the full system \eqref{rescaled_system} and compute the critical manifold $C_0$. We use the symbol ` $\cdot$ ' to denote derivatives with respect to the slow time $t$ and the symbol ` $'$ ' to denote derivatives with respect to the fast time $\tau$.

\subsubsection{Slow reduced system}\label{section_slowSubsystemAnalysis}

We obtain the slow subsystem of the full system \eqref{rescaled_system} by considering the singular limit $\eps\to 0$. In this limit, we obtain the slow reduced system
\begin{align}
  	\dot{p_1} &= (1 - q\,z)\,p_1\,, \nonumber\\
  	\dot{p_2} &= (r - (1-q)z)\,p_2\,, \nonumber\\
  	\dot{z} &= (q\,p_1 + (1-q)p_2-1)\,m\,z\,,	\label{slowSubSystem_model}\\
  	0 &= q(1-q)\,(p_1 - p_2)\,.\nonumber
\end{align}
As we describe in \eqref{slow_reduced_problem} in Appendix \ref{appendix_gspt}, this is a differential-algebraic system.


\subsubsection{Fast reduced system}\label{section_fastSubsystemAnalysis}

We scale the slow time $t$ in \eqref{rescaled_system} by $\eps$ and reformulate its dynamics in terms of the fast time $\tau=t/\eps$ to obtain
\begin{align}
	\frac{\text{d} p_1}{\text{d} \tau}&={p_1}'=\eps(1 - q\,z)\,p_1\,,\nonumber \\
	\frac{\text{d} p_2}{\text{d} \tau}&={p_2}'=\eps(r - (1-q)z)\,p_2\,, \nonumber\\
	\frac{\text{d} z}{\text{d} \tau}&={z}'=\eps(q\,p_1 + (1-q)p_2-1)\,m\,z\,,	 \label{fastSubSystem_model}\\
	\frac{\text{d} q}{\text{d} \tau}&={q}'=q(1-q)\,(p_1 - p_2)\,.\nonumber
\end{align} 
Taking the limit $\eps \to 0$ yields $p_1' = 0$, $p_2'=0$, and $z'=0$. The reduced fast system is one-dimensional, and it determines the fast dynamics of $q$ through
\begin{equation}
	{q}'=q(1-q)(p_1-p_2)\,.	\label{fast_q_dynamics}
\end{equation}
We can solve \eqref{fast_q_dynamics} explicitly for $q(\tau)$ because $p_1$ and $p_2$ are constant. For $q \in (0,1)$,
we obtain
\begin{equation}\label{equation_qh}
 	q_h(\tau;p_1,p_2) = \frac{e^{(p_1-p_2)\,\tau}}{e^{(p_1-p_2)\,\tau}+1} = \frac{1}{2}+\frac{1}{2}\tanh\left(\frac{1}{2}(p_1 - p_2) \tau\right)\,,
\end{equation}
which is gauged such that $q_h(0) = \frac{1}{2}$. 


\subsubsection{Critical manifold}\label{section_theCritManifold_fastslowSys}

The critical manifold $C_0$ is defined by the algebraic part of the slow reduced system \eqref{slowSubSystem_model}. It is given by
\begin{equation}\label{def_C0}
 	C_0 = \left\{(p_1,p_2,z,q) \in \mathbb{R}^4\,\middle\vert\, q(1-q)(p_1-p_2)=0 \right\} = \mathcal{M}_0 \cup \mathcal{M}_1 \cup \mathcal{M}_\text{sw}\,,
\end{equation}
where
\begin{align}
	 \mathcal{M}_0 &= \left\{(p_1,p_2,z,q) \in \mathbb{R}^4\,\middle\vert\, q=0\right\}\,,	\label{def_M0}\\
	 \mathcal{M}_1 &= \left\{(p_1,p_2,z,q) \in \mathbb{R}^4\,\middle\vert\, q=1\right\}\,,	\label{def_M1}\\
	 \mathcal{M}_\text{sw} &= \left\{(p_1,p_2,z,q) \in \mathbb{R}^4\,\middle\vert\, p_1=p_2\right\}\,.	\label{def_Msw}
\end{align}
We see that the critical manifold can be written as the union of a trio of three-dimensional hyperplanes.\footnote{ Because the third part, $\mathcal{M}_\text{sw}$, of the critical manifold does not play a role in the orbit construction in Section \ref{section_constructionOfThePeriodicOrbitSection}, we omit further analysis of $\mathcal{M}_\text{sw}$.} 


\subsubsection{\texorpdfstring{Slow flow on the hyperplane $\mathcal{M}_0$ \eqref{def_M0}}{Slow flow on the hyperplane M0}}
\label{section_slowflowOnCritMan_zero}

Observe that the hyperplane $\mathcal{M}_0$ \eqref{def_M0} is an invariant manifold for the slow reduced system \eqref{slowSubSystem_model}. Indeed, any initial condition with $q(0) = 0$ yields $q(t) = 0$ for all $t$ when evolved according to \eqref{slowSubSystem_model}. This allows us to study the flow of the slow reduced system \eqref{slowSubSystem_model} on the hyperplane $\mathcal{M}_0$ through the dynamical system
\begin{align}
	 \dot{p_1} &= p_1\,,\nonumber\\
	 \dot{p_2} &= (r-z)p_2\,,	\label{slowFlow_onqIsZero}\\
	 \dot{z} &= (p_2-1)\,m\,z\,.\nonumber
\end{align}
The dynamics of $p_1$ decouples from the variables $p_2$ and $z$, and the prey $p_1$ exhibits exponential growth. The predator $z$ and the prey $p_2$ form a Lotka--Volterra predator--prey system around the 
coexistence equilibrium $(p_2,z)=(1,r)$. We can thus introduce a conserved quantity on $\mathcal{M}_0$; it is given by
\begin{equation}
	H_0(p_2,z)=m\log p_2 - m\,p_2 + r \log z-z\,.	\label{LotkaVolterraEquation_p2z}
\end{equation}


\subsubsection{\texorpdfstring{Slow flow on the hyperplane $\mathcal{M}_1$ \eqref{def_M1}}{Slow flow on the hyperplane M1}}
\label{section_slowflowOnCritMan_one}

By the same reasoning as in Section~\ref{section_slowflowOnCritMan_zero}, the hyperplane $\mathcal{M}_1$ is an invariant manifold for the slow reduced system \eqref{slowSubSystem_model}. The flow of \eqref{slowSubSystem_model} on $\mathcal{M}_1$ is given by
\begin{align}
	\dot{p_1}&=(1-z)p_1\,, \nonumber \\
	\dot{p_2}&=r\,p_2\,, 		\label{slowFlow_onqIsOne} \\ \nonumber
	\dot{z}&=(p_1-1)\,m\,z\,. \nonumber \,
\end{align}
These dynamics are very similar to those on the hyperplane $\mathcal{M}_0$ \eqref{slowFlow_onqIsZero}, but the roles of $p_1$ and $p_2$ are reversed. Now the dynamics of $p_2$ decouples from the variables $p_1$ and $z$, and the prey $p_2$ exhibits exponential growth. The predator $z$ forms a Lotka--Volterra predator--prey system with the prey $p_1$ around the 
coexistence
equilibrium $(p_1,z) = (1,1)$, and the associated conserved quantity on $\mathcal{M}_1$ is given by
\begin{equation}
	H_1(p_1,z)=m\log p_1 - m\,p_1 +\log z-z\,. 		\label{LotkaVolterraEquation_p1z}
\end{equation}


\section{Construction of approximate periodic orbits}
\label{section_constructionOfThePeriodicOrbitSection}

In this section, we use the setup from Section \ref{section_preliminaries} to provide a geometric analysis of the system \eqref{rescaled_system} in terms of its slow \eqref{slowSubSystem_model} and fast \eqref{fastSubSystem_model} subsystems. We also indicate how one can construct these \textcolor{black}{singular orbits} explicitly using the analytical results in Section \ref{section_preliminaries}. \textcolor{black}{We then combine these two descriptions to construct a family of periodic orbits in the singular limit $\eps\to0$. We quantify how these singular orbits approximate solutions in the 
system \eqref{rescaled_system} for sufficiently small $\eps>0$.} In Section \ref{section_numericalContinuation}, we show numerical simulations for the constructed approximate periodic solutions for specific values of $\eps$.


\subsection{Construction of a singular periodic orbit}
\label{section_construction_singularPO}

\subsubsection{Geometric analysis of the fast reduced system}
\label{geom_analysis_fastredsys}

As we explained in Section~\ref{section_theCritManifold_fastslowSys}, the critical manifold $C_0 = \mathcal{M}_0 \cup \mathcal{M}_1 \cup \mathcal{M}_\text{sw}$ \eqref{def_C0} consists of equilibrium points of the reduced fast system \eqref{fast_q_dynamics}. Geometrically, the flow defined by \eqref{fast_q_dynamics} connects a point $(p_1,p_2,z,0) \in \mathcal{M}_0$ with a point $(p_1,p_2,z,1) \in \mathcal{M}_1$ by a heteroclinic connection, which is given explicitly by \eqref{equation_qh}. The sign of $p_1 - p_2$ determines the direction of this heteroclinic connection. Therefore, the hyperplane $\mathcal{M}_\text{sw}$ divides the four-dimensional phase space into two parts (see Figure \ref{figure:fastred}). On the side in which $p_1>p_2$, the heteroclinic connection $q_h$ \eqref{equation_qh} is directed from $\mathcal{M}_0$ to $\mathcal{M}_1$ (i.e., `upwards'); on the other side of $\mathcal{M}_\text{sw}$, in which $p_1<p_2$, the direction of the heteroclinic connection is reversed, going from $\mathcal{M}_1$ to $\mathcal{M}_0$ (i.e., `downwards'). The hyperplane $\mathcal{M}_\text{sw}$ acts as a \emph{switching plane}; when this plane is crossed, the direction of the heteroclinic flow that connects $\mathcal{M}_0$ and $\mathcal{M}_1$ is reversed.

\begin{figure}[h!]
\centering
\includegraphics[width=0.5\textwidth]{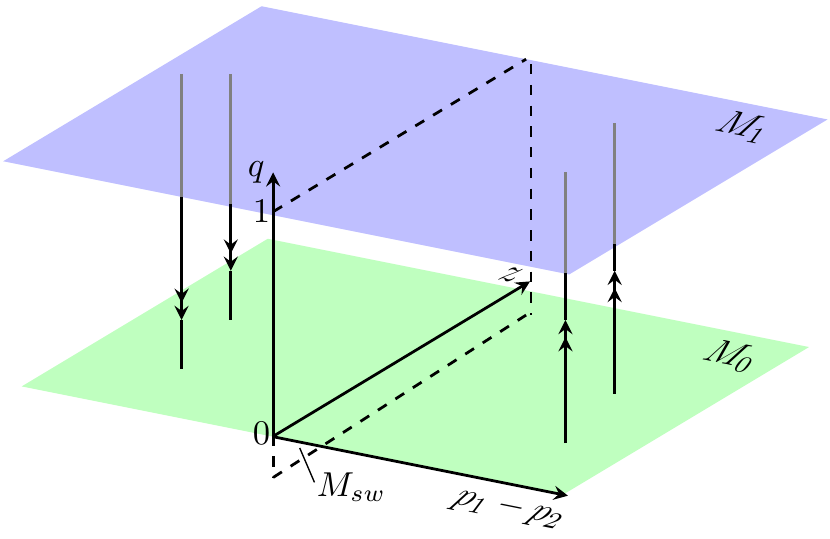} \caption{Flow of the fast reduced system \eqref{fast_q_dynamics}. For visual clarity, we depict the three slow model dimensions $(p_1,p_2,z)$ in two dimensions, which are spanned by $(p_1-p_2,z)$; the vertical axis indicates the fast variable $q$. We indicate the hyperplane $\mathcal{M}_0$ in green and the hyperplane $\mathcal{M}_1$ in blue. In this visualization, the switching hyperplane $\mathcal{M}_\text{sw}$ is spanned by the $z$-axis and $q$-axis. On the right side of $\mathcal{M}_\text{sw}$, the fast flow (indicated by double arrows) is directed upwards (i.e., from $\mathcal{M}_0$ to $\mathcal{M}_1$); on the left side of $\mathcal{M}_\text{sw}$, the direction of the fast flow is reversed (i.e., it is directed downwards).
}\label{figure:fastred}
\end{figure}


\subsubsection{Geometric analysis of the slow reduced system}\label{geom_analysis_slowredsys}

We study the flow of the slow reduced system \eqref{slowSubSystem_model} on $\mathcal{M}_0$, as given by \eqref{slowFlow_onqIsZero}, from a geometric point of view. The phase space of the flow on $\mathcal{M}_0$ is three-dimensional, and it is given in terms of the coordinates $(p_1,p_2,z)$. Projected onto the $(p_2,z)$-plane, the system \eqref{slowFlow_onqIsZero} reduces to a classical Lotka--Volterra system with conserved quantity $H_0$ \eqref{LotkaVolterraEquation_p2z}. We know \cite{Murray1993} that every orbit of this Lotka--Volterra system is closed and is determined uniquely by its value of $H_0$. Because the $p_1$ dynamics are decoupled from the $(p_2,z)$ dynamics, every $H_0$ level set in the $(p_2,z)$ plane extends to a cylindrical level set in the full $(p_1,p_2,z)$ phase space. Because $H_0$ is also a conserved quantity for the full system \eqref{slowFlow_onqIsZero}, these cylindrical level sets of $H_0$ are invariant under the flow of \eqref{slowFlow_onqIsZero}. Therefore, we can characterize the dynamics of \eqref{slowFlow_onqIsZero} by describing the three-dimensional phase space as a concentric family of cylindrical level sets of $H_0$ (see Figure \ref{figure:p1p2z_H0}).

\begin{figure}[h!]
\centering
\includegraphics[width=0.5\textwidth]{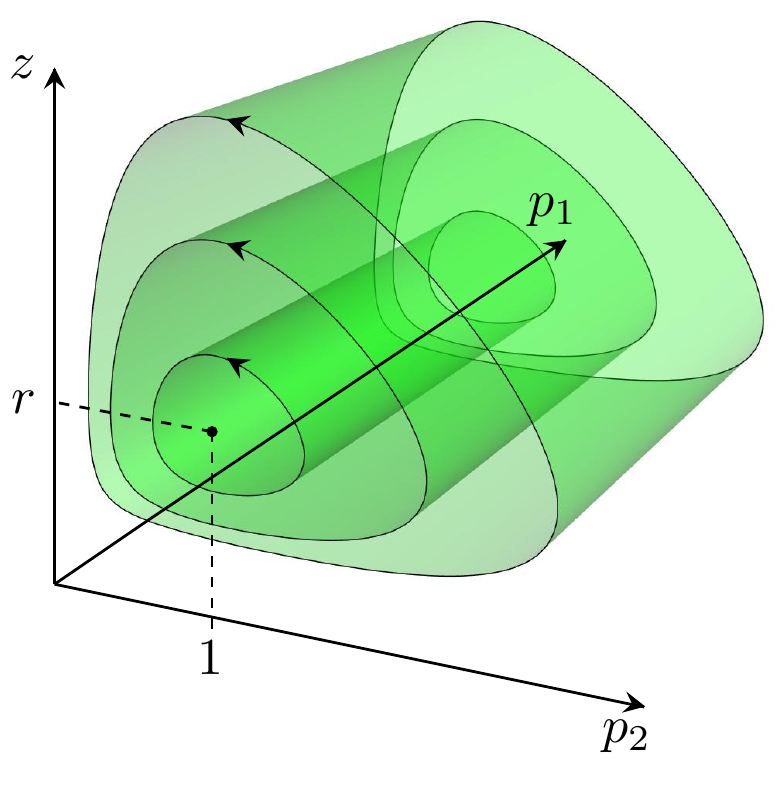}
 \caption{Flow of the slow reduced system on $\mathcal{M}_0$, as given by \eqref{slowFlow_onqIsZero}. The level sets (depicted by the green concentric cylinders) of $H_0$ \eqref{LotkaVolterraEquation_p2z} are invariant under the flow.
 }\label{figure:p1p2z_H0}
\end{figure}

For the flow of the slow reduced system \eqref{slowSubSystem_model} on $\mathcal{M}_1$, as given by \eqref{slowFlow_onqIsOne}, a geometric perspective yields an equivalent construction, with the roles of $p_1$ and $p_2$ reversed. In this case, one can characterize the dynamics of \eqref{slowFlow_onqIsOne} on the same $(p_1,p_2,z)$ phase space through another concentric family of cylindrical level sets, which are determined by the conserved quantity $H_1$ \eqref{LotkaVolterraEquation_p1z} (see Figure \ref{figure:p1p2z_H1}).

\begin{figure}[h!]
\centering
\includegraphics[width=0.5\textwidth]{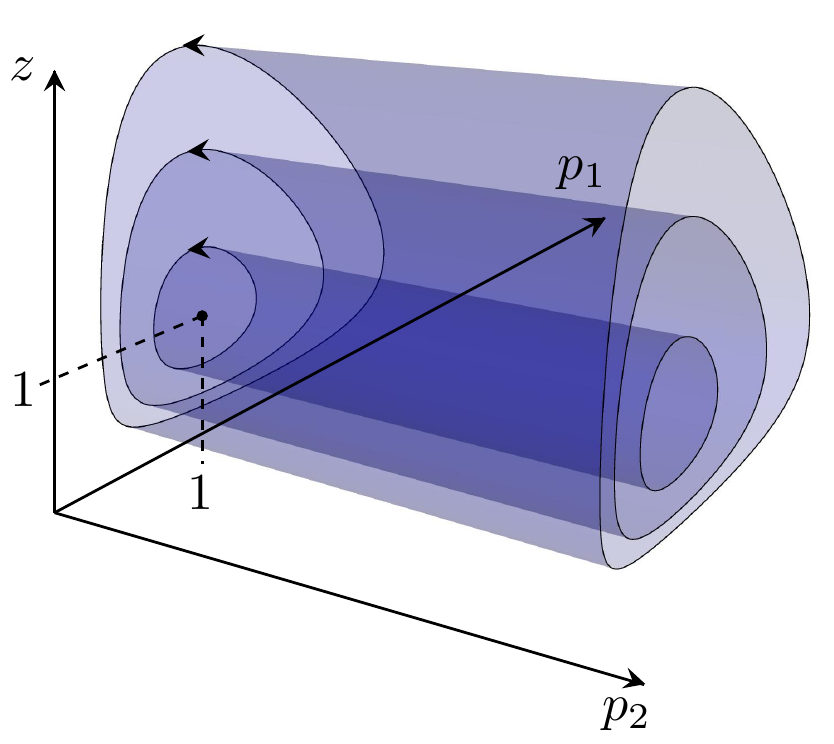}
 \caption{Flow of the slow reduced system on $\mathcal{M}_1$, as given by \eqref{slowFlow_onqIsOne}. The level sets (depicted by the blue concentric cylinders) of $H_1$ \eqref{LotkaVolterraEquation_p1z} are invariant under the flow.
  }\label{figure:p1p2z_H1}
\end{figure}


\subsubsection{Combining the fast and slow reduced dynamics}\label{section_combineslowfast}

We seek to use the geometric insights from the fast and slow reduced limits of the full system \eqref{rescaled_system} (see Sections \ref{geom_analysis_fastredsys} and \ref{geom_analysis_slowredsys}) to construct a singular periodic orbit. The idea is to exploit the heteroclinic connections between $\mathcal{M}_0$ and $\mathcal{M}_1$ on both sides of the switching plane $\mathcal{M}_\text{sw}$.

Consider a point $A_0 = (p_1^A,p_2^A,z^A,0) \in \mathcal{M}_0$, with $p_1^A>p_2^A$ (see Figure \ref{figure:per_fastred}). On this side of the switching plane $\mathcal{M}_\text{sw}$, the heteroclinic connection $q_h$ \eqref{equation_qh} takes us `up' to the corresponding point $A_1 = (p_1^A,p_2^A,z^A,1) \in \mathcal{M}_1$. We now consider the slow reduced limit, and we use the point $A_1=(p_1^A,p_2^A,z^A,1)$ as an initial condition for the slow flow $\Phi_1^t$ on $\mathcal{M}_1$. We let the slow flow on $\mathcal{M}_1$ act for some time $T_1$, so the point $A_1 \in \mathcal{M}_1$ flows to the point $B_1 = \Phi_1^{T_1} A_1 \in \mathcal{M}_1$. For notational brevity, we write $\Phi_1^{T_1} p_1^A = p_1^B$, and we use similar notation for $p_2^A$ and $z^A$, writing $B_1 = (p_1^B,p_2^B,z^B,1)$. We choose $T_1$ so that the $p_1$ coordinate is now smaller than the $p_2$ coordinate (i.e., $p_1^B < p_2^B$).

Switching back to the case of the fast reduced limit, we see (because $p_1^B < p_2^B$) that the slow flow on $\mathcal{M}_1$ has brought us to the other side of the switching plane $\mathcal{M}_\text{sw}$. We can thus use the heteroclinic connection $q_h$ to travel `down' from $B_1 \in \mathcal{M}_1$ to the corresponding point $B_0 = (p_1^B,p_2^B,z^B,0) \in\mathcal{M}_0$. Back on $\mathcal{M}_0$, we again consider the slow reduced limit, and we take the point $B_0 = (p_1^B,p_2^B,z^B,0)$ as an initial condition for the slow flow $\Phi_0^t$ on $\mathcal{M}_0$. We let the slow flow on $\mathcal{M}_0$ act for some time $T_0$. Because our goal is to construct a periodic (and hence closed) orbit, we want to choose $(p_1^A,p_2^A,z^A)$ and the times $T_0$ and $T_1$ so that the slow flow on $\mathcal{M}_0$ takes $B_0$ back to the starting point $A_0$. 

We give a schematic overview of the above construction in the following diagram: 
 \begin{equation}
  \begin{tikzcd}\label{eq:po_struc}
   \mathcal{M}_0 \arrow[twoheadrightarrow]{r}{q_h} & \mathcal{M}_1 \arrow{d}{\Phi_1^{T_1}} \\
   \arrow{u}{\Phi_0^{T_0}} \mathcal{M}_0 & \arrow[swap,twoheadrightarrow]{l}{q_h}\mathcal{M}_1
  \end{tikzcd}
\end{equation}
For a sketch of the periodic orbit in different visualizations of the phase space, see Figures \ref{figure:per_fastred}, \ref{figure:per_p1p2z}, and \ref{figure:dualphaseplane}.

\begin{figure}[h!]
\centering
\includegraphics[width=0.5\textwidth]{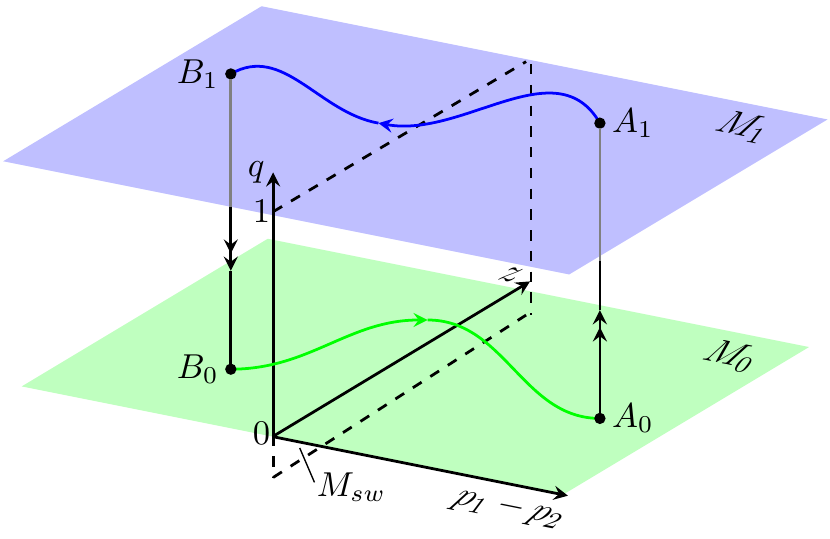}
 \caption{A visualization of the construction of the singular periodic orbit in Section \ref{section_combineslowfast} in the projection of Figures \ref{figure:p1p2z_H0} and \ref{figure:p1p2z_H1}.
 }\label{figure:per_fastred}
\end{figure}

\begin{figure}[ht!]
\centering
\includegraphics[width=0.5\textwidth]{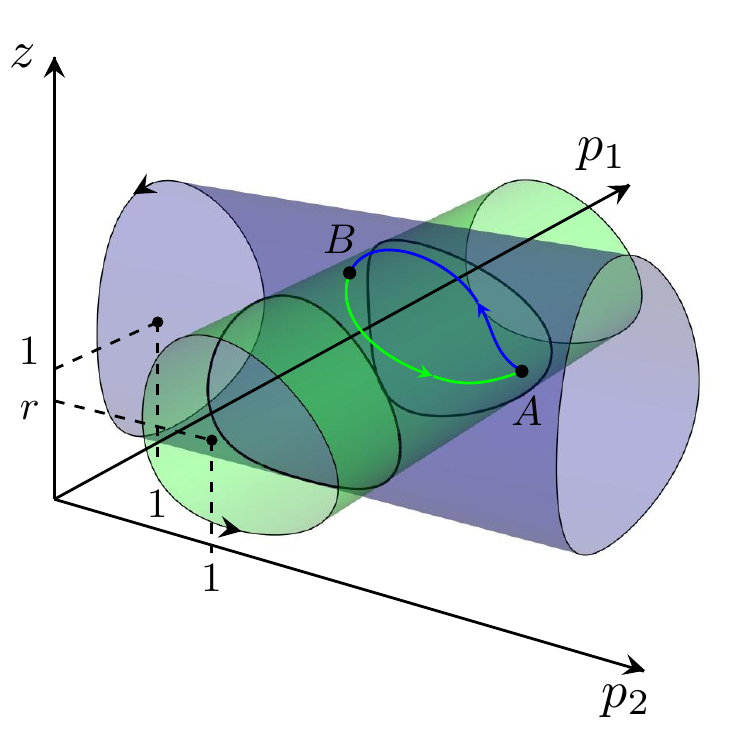}
 \caption{A visualization of the construction of the singular periodic orbit in Section \ref{section_combineslowfast} in the projection of Figure \ref{figure:fastred}.
 }\label{figure:per_p1p2z}
\end{figure}


\subsection{Existence conditions and solution families}\label{section_existenceconditions_solutionfamilies}

The goal of this section is to establish (algebraic) conditions for the existence of the closed singular orbit that we described in Section \ref{section_combineslowfast}. The periodicity of the orbit in the fast $q$ coordinate is satisfied by construction through the heteroclinic connections between $\mathcal{M}_0$ and $\mathcal{M}_1$. Therefore, we only need conditions on the initial values of the slow coordinates $(p_1,p_2,z)$ and on the slow-evolution times $T_0$ and $T_1$. We need to choose these five unknowns in a way that ensures periodicity in the three slow coordinates. Therefore, we generically expect to obtain a (possibly empty) two-parameter family of periodic orbits.\footnote{Actually, we will establish four independent conditions on the slow coordinates of the points $A_{0,1}$ and $B_{0,1}$. (In other words, there are six unknowns.)} See Figure \ref{figure:dualphaseplane} for a visualization of the singular orbit that depicts the relevant quantities that we use in the following analysis.

The expression \eqref{LotkaVolterraEquation_p2z} for the conserved quantity $H_0(p_2,z)$ of the reduced slow flow on $\mathcal{M}_0$ provides a relation between the slow coordinates of the `take-off' point $A_0$ and the slow coordinates of the `touch-down' point $B_0$. We obtain the relation from $H_0(p_2^A,z^A) = H_0(p_2^B,z^B)$, which yields
\begin{equation}\label{AB_H0condition}
	m\,\log p_2^A - m \, p_2^A + r\,\log z^A - z^A = m\,\log p_2^B - m \, p_2^B + r\,\log z^B - z^B\,.
\end{equation}
Likewise, on $\mathcal{M}_1$, we use the expression \eqref{LotkaVolterraEquation_p1z} for the conserved quantity $H_1(p_1,z)$ to obtain a relation between the slow coordinates of the touch-down point $A_1$ and the take-off point $B_1$. We obtain the relation from $H_1(p_1^A,z^A) = H_1(p_1^B,z^B)$, which yields
\begin{equation} \label{AB_H1condition}
	m\,\log p_1^A - m \, p_1^A + \log z^A - z^A = m\,\log p_1^B - m \, p_1^B + \log z^B - z^B\,.
\end{equation}

We now consider the explicit form of the slow reduced dynamics on $\mathcal{M}_1$. On $\mathcal{M}_1$, we flow the point $A_1$ to the point $B_1$ using the flow \eqref{slowFlow_onqIsOne} for a time $T_1$. The flow of the $p_2$ coordinate is decoupled from the variables $p_1$ and $z$, and it is linear, so we solve for its dynamics directly to obtain
\begin{equation}\label{T1_indirect}
	p_2^B = p_2^A e^{r\, T_1}\,.
\end{equation}
The other two slow coordinates $(p_1,z)$ interact through Lotka--Volterra dynamics. We integrate the equation for $\dot{p}_1$ in \eqref{slowFlow_onqIsOne} to yield
\begin{equation}\label{T1_direct}
	T_1 = \int_{p_1^A}^{p_1^B} \frac{1}{1-z(p_1)} \frac{\text{d}p_1}{p_1}\,,
\end{equation}
where we can obtain the expression $z(p_1)$ by invoking the conserved quantity $H_1$ and inverting the relation $H_1(p_1,z) = H_1(p_1^A,z^A)$. That is, we obtain $z(p_1)$ by solving the equation
\begin{equation}\label{p1_func_z}
	m\,\log p_1 - m \, p_1 + \log z(p_1) - z(p_1) = m\,\log p_1^A - m \, p_1^A + \log z^A - z^A\,.
\end{equation}

One can treat the slow segment on $\mathcal{M}_0$ analogously. On $\mathcal{M}_0$, we flow the point $B_0$ back to the point $A_0$ using the flow \eqref{slowFlow_onqIsZero} for a time $T_0$. The flow of the $p_1$ coordinate is decoupled from the variables $p_2$ and $z$, and it is also linear, so we can solve for its dynamics directly to obtain
\begin{equation}\label{T0_indirect}
	p_1^A = p_1^B e^{T_0}\,.
\end{equation} 
On $\mathcal{M}_0$, the slow coordinates $(p_2,z)$ interact through Lotka--Volterra dynamics. We integrate the equation for $\dot{p_2}$ using \eqref{slowFlow_onqIsZero} to yield
\begin{equation}\label{T0_direct}
	T_0 = \int_{p_2^B}^{p_2^A} \frac{1}{r-z(p_2)} \frac{\text{d}p_2}{p_2}\,,
\end{equation}
where we obtain the expression $z(p_2)$ using the conserved quantity $H_0$ by inverting the relation $H_0(p_2,z) = H_0(p_2^A,z^A)$. That is, we obtain $z(p_2)$ by solving the equation
\begin{equation}\label{p2_func_z}
	m\,\log p_2 - m \, p_2 + r\,\log z(p_2) - z(p_2) = m\,\log p_2^A - m \, p_2^A + r\,\log z^A - z^A\,.
\end{equation}

We can use the above results on the slow flow on $\mathcal{M}_0$ and $\mathcal{M}_1$ to eliminate $T_0$ and $T_1$. Combining \eqref{T1_indirect} with \eqref{T1_direct} yields
\begin{equation}\label{AB_intp1condition}
	\frac{1}{r}\log \left(\frac{p_2^B}{p_2^A}\right) = \int_{p_1^A}^{p_1^B}\frac{1}{1-z(p_1)} \frac{\text{d}p_1}{p_1}\,,
\end{equation}
and combining \eqref{T0_indirect} with \eqref{T0_direct} yields
\begin{equation}\label{AB_intp2condition}
	\log \left(\frac{p_1^A}{p_1^B}\right) = \int_{p_2^B}^{p_2^A} \frac{1}{r-z(p_2)} \frac{\text{d}p_2}{p_2}\,.
\end{equation}
Together with \eqref{AB_H0condition} and \eqref{AB_H1condition}, we now have four equations for six unknowns, which consist of the slow components of $A_{0,1}$ and $B_{0,1}$.

The relations \eqref{p1_func_z} and \eqref{p2_func_z}, which need to be inverted to obtain the integrands of \eqref{T1_direct} and \eqref{T0_direct}, are not bijective. Therefore, the specific forms of $z(p_1)$ \eqref{T1_direct} and $z(p_2)$ \eqref{T0_direct} depend on the characteristics of the underlying slow-orbit segment. To obtain computable expressions for the coordinate values of $(p_1^A,p_2^A,z^A)$ and $(p_1^B,p_2^B,z^B)$, it is necessary to characterize these underlying slow-orbit segments in more detail. We demonstrate this procedure in Appendix \ref{appendix:findsols}. We numerically evaluate the explicit expressions that we thereby obtain for the coordinate values of $(p_1^A,p_2^A,z^A)$ and $(p_1^B,p_2^B,z^B)$ for several choices of the model parameters $r$ and $m$. In the SI, we show the results of these numerical evaluations, together with visualizations of the associated singular periodic orbits.

\begin{figure}[t!]
\centering
\includegraphics[width=0.7\textwidth]{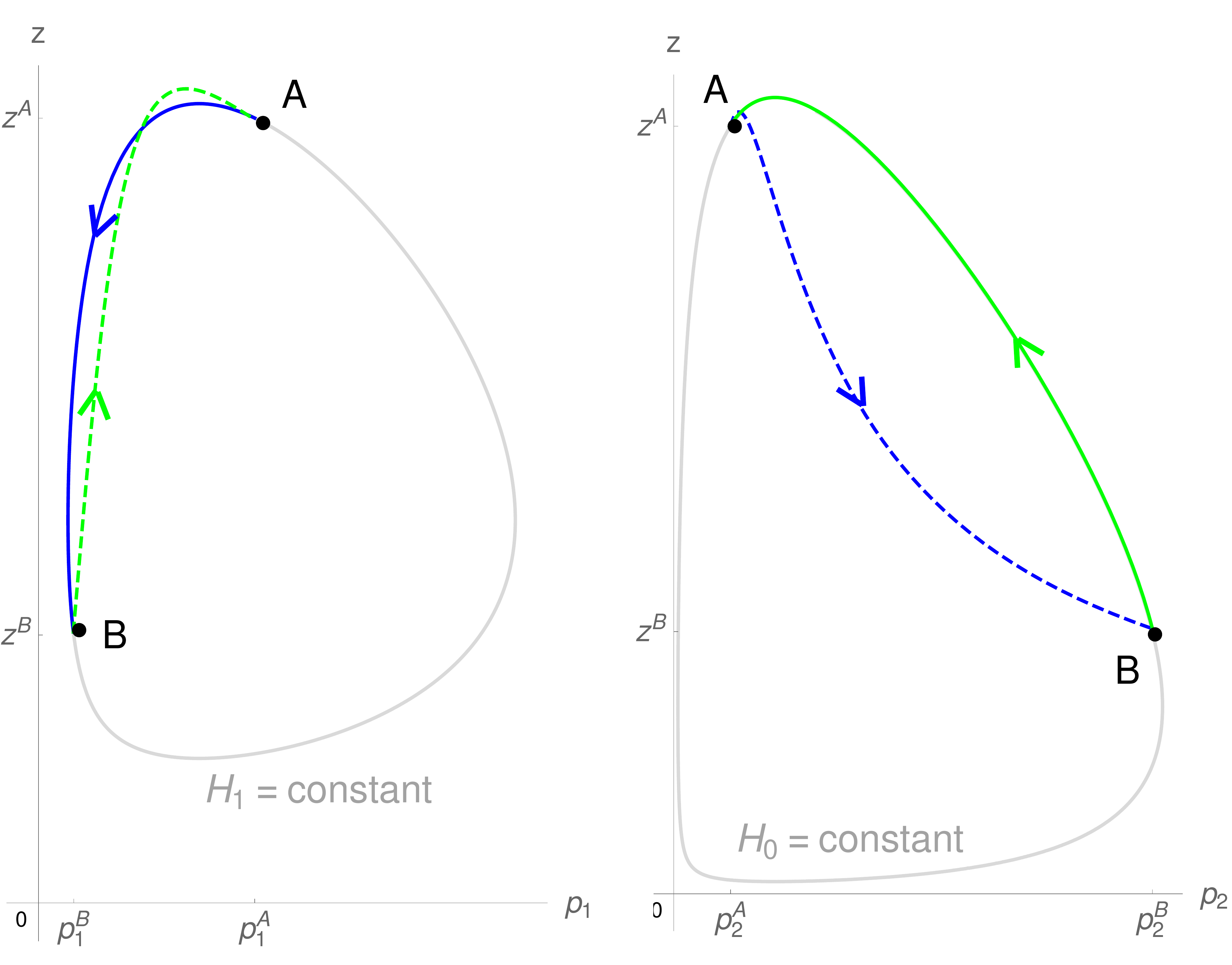}
\caption{{A dual phase-plane picture for a singular periodic orbit, as constructed in Section \ref{section_constructionOfThePeriodicOrbitSection}. We indicate the dynamics on $\mathcal{M}_1$ in blue \textcolor{black}{(solid curve on the left and dashed curve on the right)} and the dynamics on $\mathcal{M}_0$ in green \textcolor{black}{(dashed curve on the left and solid curve on the right)}. (Compare Figures \ref{figure:per_fastred} and \ref{figure:per_p1p2z}.) On the left, we show the phase plane spanned by $(p_1,z$); on the right, we show the phase plane spanned by $(p_2,z)$. In each phase plane, we use a solid curve to indicate the orbit segment in which the displayed model variables are interacting through Lotka--Volterra dynamics. We indicate the associated level curves of (left) $H_1(p_1,z)$ \eqref{LotkaVolterraEquation_p1z} and (right) $H_0(p_2,z)$ \eqref{LotkaVolterraEquation_p2z} in grey. We use a dashed curve to indicate the remaining orbit segment, in which the displayed prey variable grows exponentially. We use black dots to indicate the jump points $A$ and $B$, at which the periodic orbit `jumps' from $\mathcal{M}_0$ to $\mathcal{M}_1$, and vice versa. The arrows on the orbit segments give the direction of time.}}\label{figure:dualphaseplane}
\end{figure}


\subsection{\texorpdfstring{Approximate periodic orbits for $\eps >0$}{Approximate periodic orbits for eps > 0}}\label{section:existence}

\textcolor{black}{In Section \ref{section_construction_singularPO}, we constructed a singular periodic orbit by concatenating several orbit sections that we obtained by studying the reduced slow \eqref{slowFlow_onqIsZero}, \eqref{slowFlow_onqIsOne} and fast \eqref{fast_q_dynamics} limits of the full system \eqref{rescaled_system}. We now use these singular orbits to construct an approximation of orbits in the full system \eqref{fastslowFullSystem}.}

\begin{result}[\textbf{Existence of approximate periodic orbits}]\label{thm:existence}
 \textcolor{black}{Let $\eps > 0$ be sufficiently small, and let the coordinate triples $(p_1^A,p_2^A,z^A)$ and $(p_1^B,p_2^B,z^B)$ be such that a singular periodic orbit can be constructed according to the method outlined in Sections \ref{section_construction_singularPO} and \ref{section_existenceconditions_solutionfamilies}. Denote this singular orbit by $\gamma_0$. There then exists a solution $\gamma_\eps(t)$ of \eqref{fastslowFullSystem} and an $\mathcal{O}(1)$ time $t_*$ such that $\gamma_\eps(t)$ stays $\mathcal{O}(\eps)$ close to $\gamma_0$ for all $t \in (0,t_*)$.
}
\end{result}

\textcolor{black}{One can obtain Result \ref{thm:existence} from a mostly (though not entirely) straightforward application of `classical' perturbation theory (see, e.g., \cite{Verhulst2005}). However, because the reduced fast system \eqref{fast_q_dynamics} is one-dimensional, the slow manifolds $\mathcal{M}_0$ and $\mathcal{M}_1$ lose their locally attractive/repelling properties at the intersection with $\mathcal{M}_\text{sw}$ (for which $p_1 - p_2 = 0$), where the system exhibits a slow passage through a transcritical bifurcation \cite{KrupaSzmolyan2001}. Near $\mathcal{M}_0 \cap \mathcal{M}_\text{sw}$, we can use the standard blowup transformation $t = \sqrt{\eps} \,\tilde{t}$, $p_1 - p_2 = \sqrt{\eps} \,\tilde{p}$, $q = \sqrt{\eps} \,\tilde{q}$ to obtain, up to $\mathcal{O}(\sqrt{\eps})$, the equations 
 \begin{align*}
 	\frac{\text{d}}{\text{d} \tilde{t}}(p_1+p_2) = 0 = \frac{\text{d} z}{\text{d} \tilde{t}} 
\end{align*}
and
 \begin{align*}
  \frac{\text{d} \tilde{p}}{\text{d} \tilde{t}} &= \frac{1}{2}(1-r+z) (p_1+p_2) = \mathcal{O}(1)\,,\\
  \frac{\text{d} \tilde{q}}{\text{d} \tilde{t}} &= \tilde{p}\,\tilde{q}\,.
 \end{align*}
 It is clear that there is no exchange of stability \cite{LebovitzSchaar1975,LebovitzSchaar1977}. There is one canard, which is maximal and is given by $\tilde{q} = 0$. 
 Furthermore, using the results in \cite{DeMaesschalckSchecter2015}, it follows that $0<\tilde{q}(\tilde{t})<\tilde{q}(-a)$ for all $-a<\tilde{t}<a$ with $a \in \mathcal{O}(1)$. In other words, $\tilde{q}(\tilde{t})$ stays close to $\tilde{q} = 0$ for $\mathcal{O}(1)$ time in $\tilde{t}$. The analysis at $\mathcal{M}_1 \cap \mathcal{M}_\text{sw}$ is analogous; for more details on the local analysis near $p_1 - p_2 = 0$, see \cite{DeMaesschalckSchecter2015}. One can apply the classical theory \cite{Verhulst2005}, which guarantees the existence of a solution to \eqref{fastslowFullSystem} that is $\mathcal{O}(\eps)$ close to the singular approximation $\gamma_0$, on either side of $\mathcal{M}_\text{sw}$. The above blowup argument shows that these classical solutions stay close to either $\mathcal{M}_{0}$ or $\mathcal{M}_{1}$ when crossing $\mathcal{M}_\text{sw}$.}

\begin{rmk}\label{rmk:noSotoTrevino}
 \textcolor{black}{In the context of geometric singular perturbation theory, the standard reference for the existence of periodic orbits constructed by concatenating slow and fast orbit segments (as outlined in Section \ref{section_construction_singularPO}) is the seminal paper by Soto-Trevi{\~n}o \cite{Soto-Trevino2001}. However, the system \eqref{fastslowFullSystem} that we analyze in our current paper has only one fast component. Consequently, one cannot use the standard notion of `normal hyperbolicity,' because its definition requires the number of normal directions to be at least two. This, in turn, implies that one cannot apply the theory from \cite{Soto-Trevino2001} to the case at hand. Moreover, the existence of a two-parameter family of singular periodic orbits indicates that the intersection of the stable and unstable manifolds of $\mathcal{M}_{0,1}$ is not transversal in the singular limit $\eps \to 0$. Therefore, generically, a singular orbit $\gamma_0$ constructed as outlined in Section \ref{section_construction_singularPO} does \emph{not} perturb to a periodic orbit in the full system \eqref{fastslowFullSystem}. To find a proper transversal intersection of the stable and unstable manifolds of $\mathcal{M}_{0,1}$, one would need to extend the leading-order analysis presented in the present paper to higher orders in $\eps$ to obtain additional existence conditions to match the number of free parameters (see Section \ref{section_existenceconditions_solutionfamilies}).}
\end{rmk}

}


\section{Ecologically relevant qualitative aspects of the constructed periodic orbits}\label{section_QualAsp}

In this section, we discuss several qualitative aspects of periodic solutions that are ecologically relevant---including synchronization between predator and prey and/or between two prey species, clockwise cycles, and counterclockwise cycles. \textcolor{black}{For a summary of cyclic behavior exhibited by the singular periodic orbits constructed from the model in \eqref{fastslowFullSystem}, see Table \ref{summary_cyclicBehavior}.} These types of behavior occur (1) in the data collected from microscopic aquatic organisms in field research \cite{TirokGaedke2006,TirokGaedke2007a} and under laboratory conditions \cite{yoshida2007cryptic,Becks2010,HiltunenBecksNatComm2014} and (2) in experimental studies of coevolution in phage--bacteria systems \cite{Mizoguchi2003,Weietal2010}. Importantly, such behavior also arises in the family of periodic orbits that we constructed in Section \ref{section_constructionOfThePeriodicOrbitSection}. 

\begin{table}[ht!]
\centering
\textcolor{black}{
 \begin{tabular}{p{6cm} p{7cm} r}
 Name & Description & Figure number \\ \hline
 Prey--prey synchronization & The two prey oscillate in antiphase & \ref{figure:appendix_fig_antiphase}\\ \hline
 Predator--prey--prey synchronization & The predator alternates between (1) being in phase with prey 1 and in antiphase with prey 2 and (2) being in antiphase with prey 1 and in phase with prey 2 & \ref{figure:sync_predpreyprey}\\ \hline
 Predator--prey synchronization & The predator alternates between (1) being in phase with prey 2 and (2) being in antiphase with prey 2 & \ref{figure:sync_predp2}\\ \hline
 Counterclockwise cycles & Prey peaks before predator & \ref{figure:sync_predpreyprey}, \ref{figure:sync_predp2}\\ \hline
 Clockwise cycles & Predator peaks before prey & \ref{figure:appendix_fig_clockwise}, \ref{figure:appendix_fig_hybridphase_mid}
  \end{tabular}
\caption{\textcolor{black}{Summary of the possible oscillatory behavior exhibited by the singular periodic orbits that we construct for the 1 fast--3 slow system in \eqref{rescaled_system}.}}
\label{summary_cyclicBehavior}
}
\end{table}

We have access to field data on microscopic aquatic organisms \cite{TirokGaedke2006,TirokGaedke2007a}, and these data exhibit both antiphase and in-phase oscillations between the two different prey types. \textcolor{black}{In terms of the fast--slow 1 predator--2 prey system \eqref{fastslowFullSystem}, these data can be used to obtain values for the model parameters---in particular, the prey growth rates $r_1$ and $r_2$ (which determine $r$ in the rescaled system in \eqref{rescaled_system}), the predator conversion efficiency ($e$), and the predator death rate ($m$). All values chosen in the current paper for these model parameters are within the range suggested in previous modeling work that used these data (see, e.g., \cite{TirokGaedke2010,ourSmoothPaper}).} \textcolor{black}{We also take into account experimental evidence of prey preference exhibited by the predator species in these data \cite{MullerSchlegel1999} and assume that predator $z$ prefers $p_1$ and thus can exert more grazing pressure on it than on its alternative prey $p_2$. This difference in prey preference manifests in the model in \eqref{fastslowFullSystem} via the parameter $q_2$, which scales the benefit that the predator obtains from feeding on $p_2$. Such an advantage of experiencing lower predation pressure can be explained, for example, by a difference in the use of limited nutrients between the two different prey.} \textcolor{black}{The alternative prey could, for example, invest resources in building defense mechanisms (such as a hard silicate cover) that make it difficult for the predator to digest this alternative prey.} \textcolor{black}{In this example, the preferred prey $p_1$ has a poorer defense than the alternative prey $p_2$; as a trade-off, prey $p_1$ has more resources available than $p_2$ to be used for population growth, as these resources are not invested in building defense mechanisms.}
\textcolor{black}{Consequently, we follow earlier modeling work on these data \cite{TirokGaedke2010} and incorporate such a prey \emph{trade-off} in the model by assuming that the growth rate of the preferred prey is larger than that of the alternative prey. That is, we assume that $r_1>r_2$ in \eqref{fastslowFullSystem}, which implies that $0<r<1$ in the rescaled system in \eqref{rescaled_system}.}

Although there are several experimental studies of coevolution in microscopic aquatic organisms \cite{Fussmann2000,Yoshidaetal2003,HiltunenBecksNatComm2014} (and in phage--bacteria systems \cite{Weietal2010,Mizoguchi2003,Halletal2011}, which is another type of exploiter--resource system that translates to a model of predator--prey interaction), we have not yet encountered empirical observations of evolutionary and demographic dynamics in a system of one predator and two different prey species. 
In the bacterium--phage system studied in \cite{Weietal2011}, a bacterial (i.e., prey) subpopulation that replicates slowly and is phenotypically resistant to the phage (i.e., predator) was suggested as a possible mechanistic explanation for the observed dynamics. This is an example of our model assumption of the growth rate of the alternative prey being less than that of the preferred prey (i.e., $0<r<1$ in the rescaled system in \eqref{rescaled_system}).


\subsection{Synchronization}\label{synchronizationSubSection}


\subsubsection{Prey--prey synchronization}\label{section:preypreyphase}

When two populations oscillate in phase or in antiphase, local extrema of two species occur at exactly the same instances in time. Because $p_1$ increases monotonically on $\mathcal{M}_0$ and $p_2$ increases monotonically on $\mathcal{M}_1$, it follows that $p_1$ does not have an extremum during the slow dynamics on $\mathcal{M}_0$ and that $p_2$ does not have an extremum during the slow dynamics on $\mathcal{M}_1$. Therefore, for $p_1$ and $p_2$ to oscillate in phase or in antiphase, we need their local, aligned extrema to occur at the jump points---i.e., where $q$ jumps from $0$ to $1$, or vice versa. Because $p_1$ increases from jump point $B$ towards jump point $A$, it follows that $p_1$ has a local maximum at $A$ and a local minimum at $B$. By the same reasoning, $p_2$ must have a local minimum at $A$ and a local maximum at $B$. Therefore, we can conclude that 
\emph{the only type of prey--prey synchronization that occurs in the singular periodic solutions that we have constructed is when the two prey species oscillate in antiphase.} Therefore, based on the above considerations, the dual phase-plane picture associated with prey--prey synchronization must be as depicted in Figure \ref{figure:appendix_fig_antiphase}. As is clear from this figure, the $z$ coordinates of both jump points $A$ and $B$ must lie above the $z=1$ nullcline. That is, $z_A > 1$ and $z^B > 1$. Our numerical calculations show that there exists a two-parameter family of periodic orbits in which both prey species oscillate in antiphase for a range of the model parameters $(r,m)$; see Figure \ref{figure:exconds_sample} in Appendix \ref{appendix:findsols} for a visualization of such a solution family. See the SI for visualizations of the associated singular periodic orbits and for our numerically obtained values of $(p_1^A,p_2^A,z^A)$ and $(p_1^B,p_2^B,z^B)$.

\begin{figure}[t!]
\centering
\includegraphics[width=0.7\textwidth]{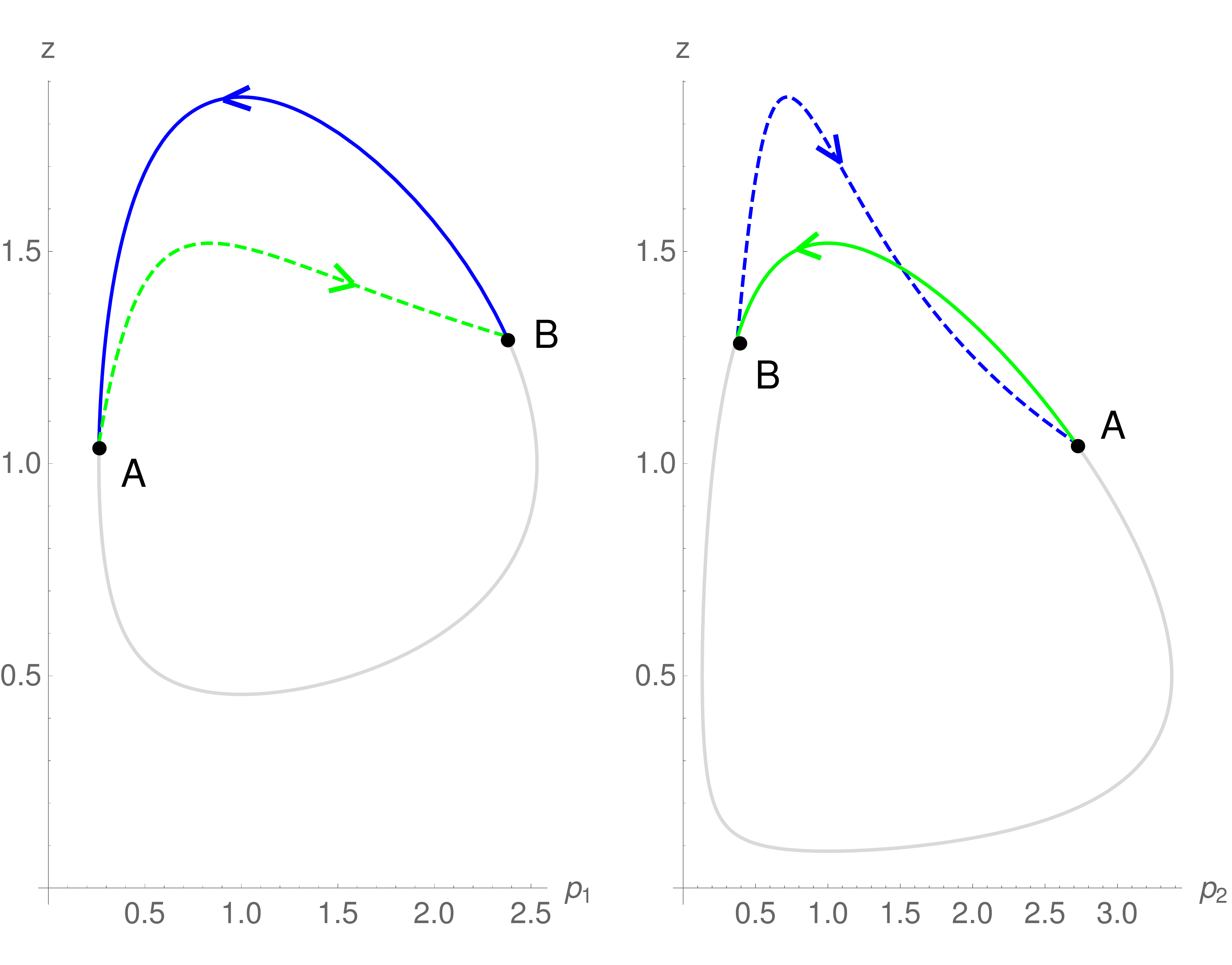}
\caption{The dual phase-plane picture for a singular periodic orbit with prey--prey synchronization. In this example, $(r,m) = (0.5,0.4)$. The local extrema of $p_1$ and $p_2$ are located at the jump points $A$ and $B$. We indicate the dynamics on $\mathcal{M}_0$ in green \textcolor{black}{(dashed curve on the left and solid curve on the right)} and the dynamics on $\mathcal{M}_1$ in blue \textcolor{black}{(solid curve on the left and dashed curve on the right)}. (Compare Figures \ref{figure:per_fastred} and \ref{figure:per_p1p2z}.)}\label{figure:appendix_fig_antiphase}
\end{figure}


\subsubsection{Predator--prey synchronization}

Because our {model includes two prey species, predator--prey synchronization can potentially arise either through synchronization between the predator and prey 1 or though synchronization between the predator and prey 2. Predator and prey densities that oscillate almost exactly out of phase with each other have been observed in experimental studies on the effects of rapid prey evolution on ecological dynamics \cite{Becks2010,Yoshidaetal2003}. In this section, we examine the conditions under which there is a jump from one slow manifold to the other, and we thereby show that our model exhibits oscillations in which the predator is in phase with a prey at one jump point and out of phase with the same prey at the other jump point.


\paragraph{Synchronization between all three species.} 

The predator cannot oscillate either in phase or in antiphase with both prey species simultaneously, because this would imply that the aligned local extrema of the prey species are the same type (i.e., both maxima or both minima). In other words, the two prey species would exhibit in-phase oscillations, and we showed in Section \ref{section:preypreyphase} that this cannot occur. However, it is still possible for the predator to oscillate in phase with one prey and in antiphase with the other prey. Suppose that the two prey oscillate in antiphase, as described in Section \ref{section:preypreyphase}. If the predator is in phase with one prey and in antiphase with the other, then the predator density has local extrema located at the jump points. The nature of these extrema is dictated by the `jump conditions' on $p_{1,2}$ at the jump points $A$ and $B$ in the following manner. Suppose that the predator density $z$ has a local maximum at $A$. From the slow flow on $\mathcal{M}_0$ \eqref{slowFlow_onqIsZero} and the slow flow on $\mathcal{M}_1$ \eqref{slowFlow_onqIsOne}, we infer that this implies that $(p_2^A - 1)m z^A > 0$ and $(p_1^A-1)m z^A < 0$. However, this violates the jump condition that $p_1^A>p_2^A$ (see Section \ref{section_combineslowfast}). Therefore, $z$ cannot have a local maximum at $A$. By an analogous argument, we see that $z$ cannot have a local maximum at $B$. Therefore, if the predator density has local extrema at both jump points $A$ and $B$, then both of these extrema are local minima. We conclude that the only way in which the predator is synchronized with both prey species is if the predator is alternating between (1) being in phase with prey 1 and in antiphase with prey 2 and (2) being in antiphase with prey 1 and in phase with prey 2. For an example of such predator--prey--prey synchronization, see Figure \ref{figure:sync_predpreyprey}.

\begin{figure}
\center
  \includegraphics[width=0.6\textwidth]{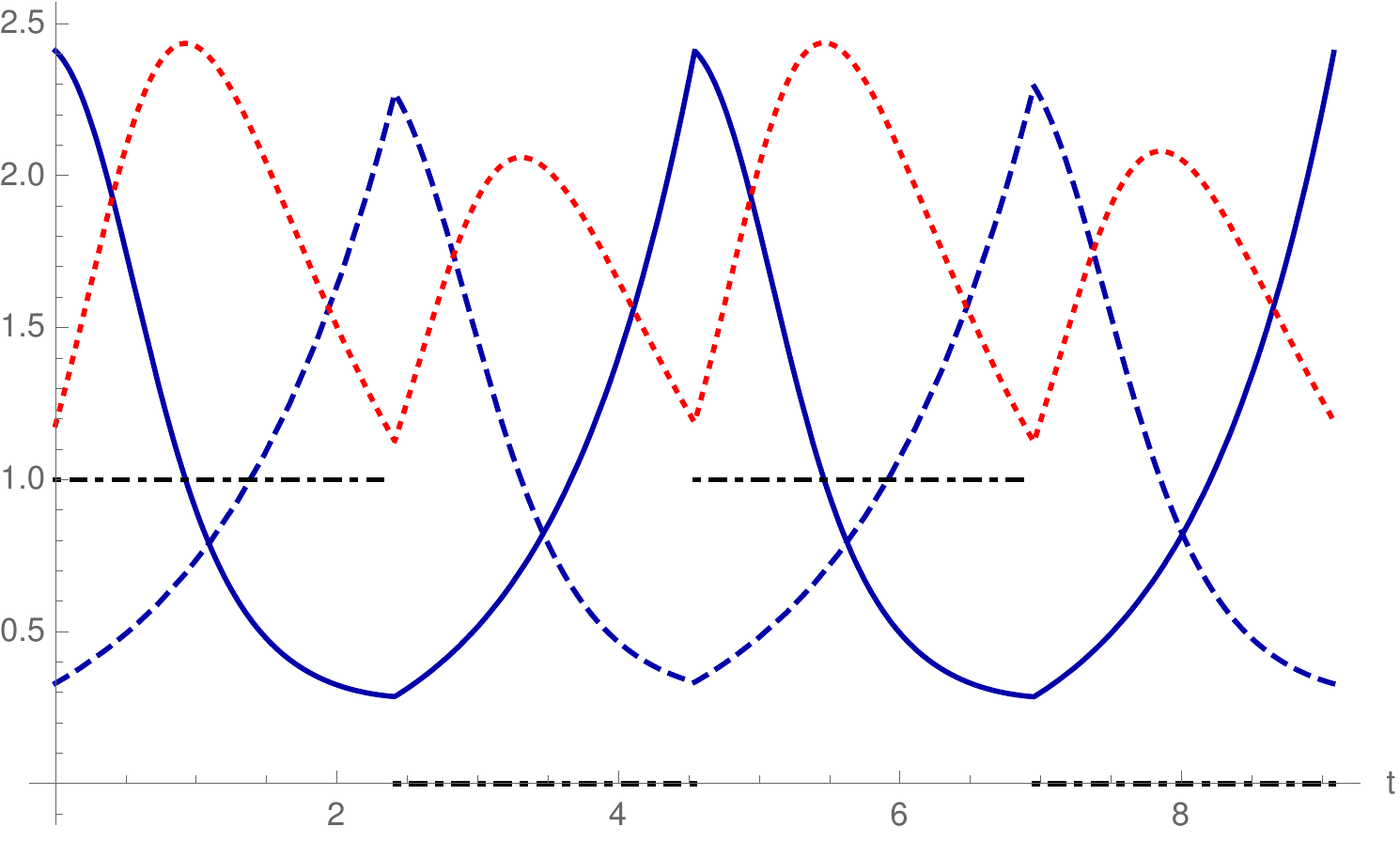}
\caption{Rescaled abundances of the preferred prey $p_1$ (solid blue curve), alternative prey $p_2$ (dashed blue curve), the predator $z$ (dotted red curve), and predator trait $q$ (dot-dashed black lines) as a function of the rescaled time $t$ (on the  horizontal axis) for a singular orbit exhibiting predator--prey--prey synchronization in the system in \eqref{rescaled_system} with $r=0.8$, $m=1$, and $\eps=0$. The jump points are located at $(p_1^A,p_2^A,z^A) \approx (2.41,0.33,1.18)$ and $(p_1^B,p_2^B,z^B) \approx (0.29,2.27,1.39)$.}
\label{figure:sync_predpreyprey}
\end{figure}


\paragraph{Synchronization between predator and one prey.} 

Suppose that the predator $z$ oscillates in phase (or in antiphase) with one prey only, which we assume is $p_1$ without loss of generality. From the slow dynamics of $p_1$ on $\mathcal{M}_1$ \eqref{slowFlow_onqIsOne} and on $\mathcal{M}_0$ \eqref{slowFlow_onqIsZero}, it is clear that during the slow dynamics, the local extrema of $z$ and $p_1$ are unable to align. On $\mathcal{M}_0$, the prey $p_1$ changes monotonically; on $\mathcal{M}_1$, the prey $p_1$ and predator $z$ are related through Lotka--Volterra dynamics, which forbids alignment of local extrema of the participating predator and prey. We thus conclude that the local extrema of the predator $z$ and the prey species $p_1$ with which it oscillates in phase (or in antiphase) must occur at the jump points $A$ and $B$. From the fact that $p_1$ increases monotonically during the slow dynamics on $\mathcal{M}_0$, we conclude that $p_1$ must be a maximum if it has a local extremum at $A$. This, in turn, implies that the derivative of $p_1$ on $\mathcal{M}_1$ at the jump point $A$ must be negative. Therefore, $(1-z^A)p_1^A < 0$ (see \eqref{slowFlow_onqIsOne}), so $z^A > 1$. For the other prey species $p_2$, we see using \eqref{slowFlow_onqIsOne} that the derivative of $p_2$ on $\mathcal{M}_1$ at $A$ is (trivially) positive and that the derivative of $p_2$ on $\mathcal{M}_0$ at $A$ is given by $(r-z^A)p_2^A$. However, we have already concluded  that $z^A>1$, and because we have also assumed that $0<r<1$ (see Section \ref{section_rescalingthesystem}), the derivative of $p_2$ on $\mathcal{M}_0$ at $A$ must be negative. Therefore, if $p_1$ has a local extremum at $A$ (which must be a maximum), then $p_2$ has a local minimum at $A$. This violates the assumption that the predator $z$ oscillates in phase (or in antiphase) with one prey only. We therefore conclude that the predator $z$ cannot oscillate in phase (or in antiphase) with prey 1 only. Any synchronization between the predator and prey 1 necessarily implies synchronization between the predator and prey 2.  We discussed this predator--prey--prey synchronization in the previous paragraph. It is worthwhile to note that the simultaneous synchronization of both prey is a direct consequence of the asymmetry between the dynamics of $p_1$ and $p_2$, manifested through the parameter $r \in (0,1)$ in \eqref{slowFlow_onqIsOne} and \eqref{slowFlow_onqIsZero}.\phantomsection\label{sentence:1vsr}

Now suppose that the predator $z$ oscillates in phase (or in antiphase) with one prey only, and suppose that that prey is $p_2$. By the same arguments as above, we can conclude that if $p_2$ has local extrema at $A$ and $B$, then $p_2$ must have a local minimum at $A$ and a local maximum at $B$. Therefore, from \eqref{slowFlow_onqIsOne} and \eqref{slowFlow_onqIsZero}, it follows that $r< z^A$ and $r < z^B$. Because we assumed that $z$ does not oscillate in phase (or in antiphase) with prey $p_1$,  the derivative of $p_1$ must be positive on both sides of jump point $A$, and it must also be positive on both sides of jump point $B$. This implies that $z^A<1$ and $z^B < 1$. We can conclude that it is possible for the predator $z$ to synchronize with one prey only, and this prey must be prey 2. This situation occurs if and only if $r<z^A<1$ and $r<z^B<1$. However, from the previous paragraph, we know that the local extrema of the predator $z$ at the jump points $A$ and $B$ can only be minima. Therefore, the predator $z$ and prey $p_2$ cannot oscillate in phase or in antiphase. The only synchronization possible between $z$ and $p_2$ is of an `alternating' type: if $z$ and $p_2$ are synchronized, then they are in phase at jump point $A$ and in antiphase at jump point $B$. For an example of this type of synchronization, see Figure \ref{figure:sync_predp2}.

\begin{figure}
\center
  \includegraphics[width=0.6\textwidth]{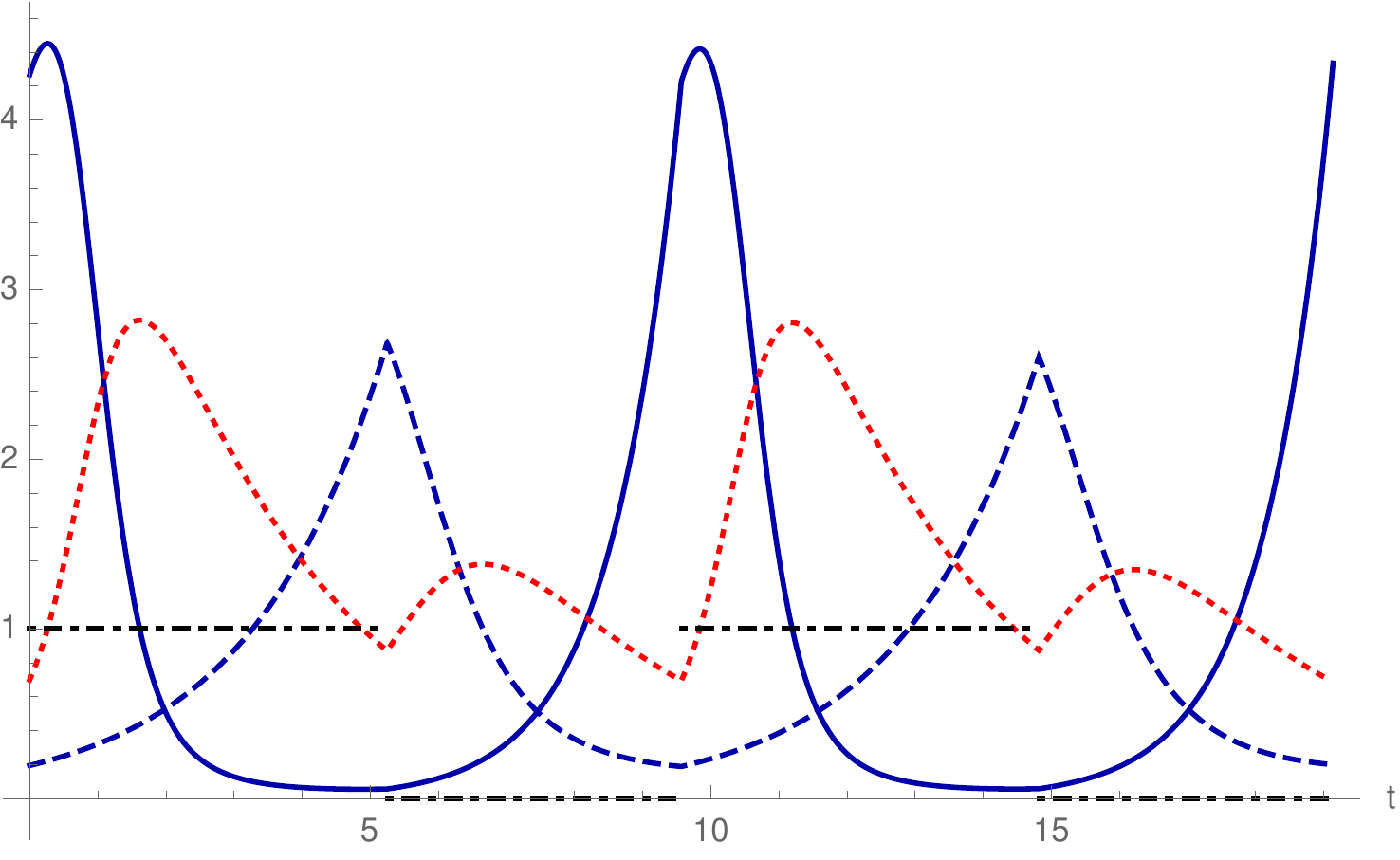}
\caption{Rescaled abundances of the preferred prey $p_1$ (solid blue curve), alternative prey $p_2$ (dashed blue curve), the predator $z$ (dotted red curve), and predator trait $q$ (dot-dashed black lines) as a function of the rescaled time $t$ (on the horizontal axis) for a singular orbit exhibiting synchronization between predator and prey 2 in the system in \eqref{rescaled_system} with $r=0.5$, $m=0.4$, and $\eps=0$. The jump points are located at $(p_1^A,p_2^A,z^A) \approx (4.27,0.19,0.7)$ and $(p_1^B,p_2^B,z^B) \approx (0.06,2.69,0.85)$. \textcolor{black}{Note that the local maxima for $p_1$ occur just after the jump point, where $p_2$ and $z$ have a synchronous local minimum. On closer inspection, one can see a sudden change in the slope of $p_1$ at that jump point. Compare this figure with Figure \ref{figure:sync_predpreyprey}, in which all species are synchronized.}
}
\label{figure:sync_predp2}
\end{figure}

Based on the above analysis, we see that the only phase-locking mode available for the periodic orbits that we have constructed in this paper is a `hybrid' phase in which the predator is in phase with one prey at one jump point and out of phase with the same prey at the other jump point. In summary, the above analysis on synchronization types yields the following insights:
\begin{enumerate}
 \item{Alignment of local extrema of the model species $(p_1,p_2,z)$ can occur only at a jump point.}
 \item{If the predator $z$ has a local extremum at a jump point, then this extremum must be a local minimum.}
 \item{When $z^A>r$ and $z^B>r$, the predator $z$ and prey $p_2$ are in phase at $A$ and in antiphase at $B$.} 
 \item{When $z^A>1$ and $z^B>1$, the predator $z$ and prey $p_1$ are in antiphase at $A$ and in phase at $B$.} 
 \item{Therefore, the only two types of synchronization between the predator $z$ and prey $p_{1,2}$ are the ones in Figures \ref{figure:sync_predpreyprey} and \ref{figure:sync_predp2}.}
\end{enumerate}


\subsection{Clockwise and counterclockwise cycles}

In this section, we discuss the ordering of the peak abundances of the predator and prey populations in cycles exhibited by the model in \eqref{rescaled_system}. In particular, we describe two types of situations: (1) a peak in the predator abundance precedes that in the prey population (so the cycles have a `clockwise' orientation when depicted on a predator--prey phase plane), and (2) a peak in the prey abundance is followed by that in the predator population (so the flow travels `counterclockwise' in that phase plane).

The nomenclature `clockwise' and `counterclockwise' stems from the orientation of the flow in a classical Lotka--Volterra system, which describes one prey species and one predator species. In the traditional phase-plane depiction of the Lotka--Volterra system, the prey is placed on the horizontal axis, and the predator is placed on the vertical axis. In this case, the Lotka--Volterra flow has a counterclockwise orientation. In the solution time series, this dynamical behavior is characterized by the fact that a peak in the prey population is relatively close \textcolor{black}{(i.e., a quarter of a period with a small perturbation from the equilibrium point)} to a peak in the predator population. Moreover, in Lotka--Volterra dynamics, the prey peaks first and the predator peaks shortly thereafter.

The prediction of counterclockwise cycles due to density-dependent predator--prey interactions in the Lotka--Volterra model is supported by empirical observations collected from hare and lynx populations \cite{eltonNicholson1942lynxhare}. As is also the case for several other traditional predator--prey models, the Lotka--Volterra system assumes that the behavior and characteristics of the organisms remain fixed on the time scale of ecological interactions. As we discussed in Section \ref{section_introduction} and will discuss further in Section \ref{synchronizationSubSection}, rapid evolution alters the population dynamics and, in particular, it can generate cycles in which the peak in the predator abundance follows the peak in the prey population with a phase lag that is larger than a quarter of a period \cite{Becks2010}.

In contrast to the counterclockwise cycles, clockwise cycles are characterized by a negative phase lag between the peak abundances and a reversed ordering of the predator and prey maxima. Recently, Cortez and Weitz \cite{CortezWeitz2014} analyzed ecological data sets collected from various predator--prey systems and identified regions in them that have a clockwise orientation. In \cite{CortezWeitz2014}, a peak in the predator population is construed to precede a peak in the prey abundance if the time between a predator peak and the following prey peak is less than the distance between the predator peak and the preceding prey peak. Modeling suggests that evolutionary changes on a time scale comparable to that of the ecological interactions and occurring in both predator and prey offer a possible mechanistic explanation for the reversed ordering of the peak abundances \cite{CortezWeitz2014}. In other words, a small population of a prey type that has invested in predator defense mechanisms and a large population of a predator that is ineffective in counteracting the prey's defense can yield low predator abundance and high prey abundance because of the effective prey defense. Consequently, selection favors predators that are effective in counteracting prey defense, so that the prey population starts to decrease. Simultaneously, the predator population remains low because of the high cost of counteracting prey defense. However, due to low predator population, there is room for the prey population with low predator defense to increase. The predator population then increases because of a high abundance of more vulnerable prey. In this reversed situation, the predator peaks first, and the prey peaks shortly thereafter. 

In the predator--prey-prey system that we study in the present paper, the answer to the question of whether clockwise cycles occur seems to be straightforward. Using the orbit visualization in Figure \ref{figure:appendix_fig_antiphase}, we immediately see that both the interaction between the predator and prey 1 and the interaction between the predator and prey 2 occur in a counterclockwise fashion. Such orientations are inherent to the construction of the (singular) periodic orbits in question (see also Figure \ref{figure:per_p1p2z}). However, from such phase portraits, one is unable to draw any conclusion about the difference in time between predator and prey extrema. In particular, using a phase-space perspective, one cannot readily deduce whether a prey peak is shortly followed by a predator peak (the time-series hallmark of a counterclockwise cycle), or vice versa.

To obtain more insight on the relative time difference between predator and prey peaks, we use the analysis of Section \ref{synchronizationSubSection}. The slow dynamics on $\mathcal{M}_0$ [see \eqref{slowFlow_onqIsZero}] and $\mathcal{M}_1$ [see \eqref{slowFlow_onqIsOne}] show that on either slow manifold, one prey species increases monotonically and the other prey species interacts with the predator through Lotka--Volterra dynamics. In Lotka--Volterra dynamics, a peak in prey density always precedes a predator peak. Therefore, the situation in which a predator peak is shortly followed by a prey peak cannot occur during the slow dynamics on either $\mathcal{M}_0$ or $\mathcal{M}_1$. However, as we saw in Section \ref{synchronizationSubSection}, local prey maxima can also occur at the jump points $A$ and $B$. (The predator can only have local minima at the jump points; see Subsection \ref{section:preypreyphase}.) In Figure \ref{figure:appendix_fig_clockwise}, we show an example of a singular periodic orbit in which the prey density $p_1$ has a local maximum at $A$, where $q$ jumps from $0$ to $1$. The predator peak occurs shortly before this instance in time. In Figure \ref{figure:appendix_fig_hybridphase_mid}, we show another singular periodic orbit, in which the predator peak occurs almost exactly in between the peak of prey 1 and the peak of prey 2. We can therefore conclude that the dynamics of our model \eqref{fastslowFullSystem} admits (singular) periodic solutions whose time series exhibit an ordering of predator and prey peaks that can be interpreted as `clockwise' in the sense that a prey peak is shortly preceded by a predator peak. In Section \ref{section_numericalContinuation}, we will show that the localization in time of these local maxima persists when we increase the value of the small parameter $\eps$.

\begin{figure}
\center
  \includegraphics[width=0.6\textwidth]{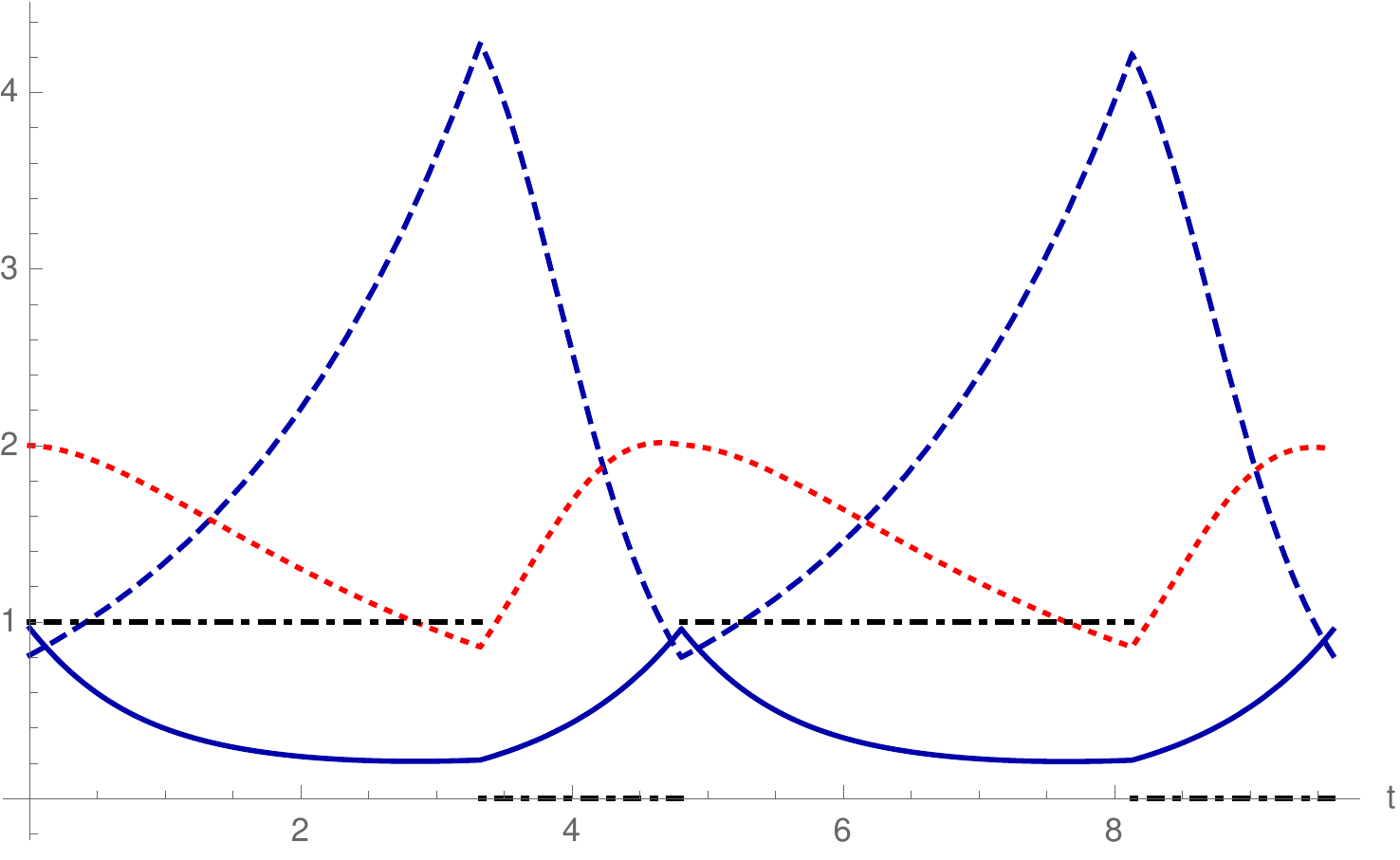}
\caption{Rescaled abundances of the preferred prey $p_1$ (solid blue curve), alternative prey $p_2$ (dashed blue curve), the predator $z$ (dotted red curve), and predator trait $q$ (dot-dashed black lines) as a function of the rescaled time $t$ (on the  horizontal axis) for a singular orbit exhibiting `clockwise' behavior (i.e., the peak in the predator density occurs just before the peak in prey 1) in the system in \eqref{rescaled_system} with $r=0.5$, $m=0.4$, and $\eps=0$. The jump points are located at $(p_1^A,p_2^A,z^A) \approx (0.97,0.81,2.0)$ and $(p_1^B,p_2^B,z^B) \approx (0.22,4.28,0.85)$.
}\label{figure:appendix_fig_clockwise}
\end{figure}

\begin{figure}
\center
  \includegraphics[width=0.6\textwidth]{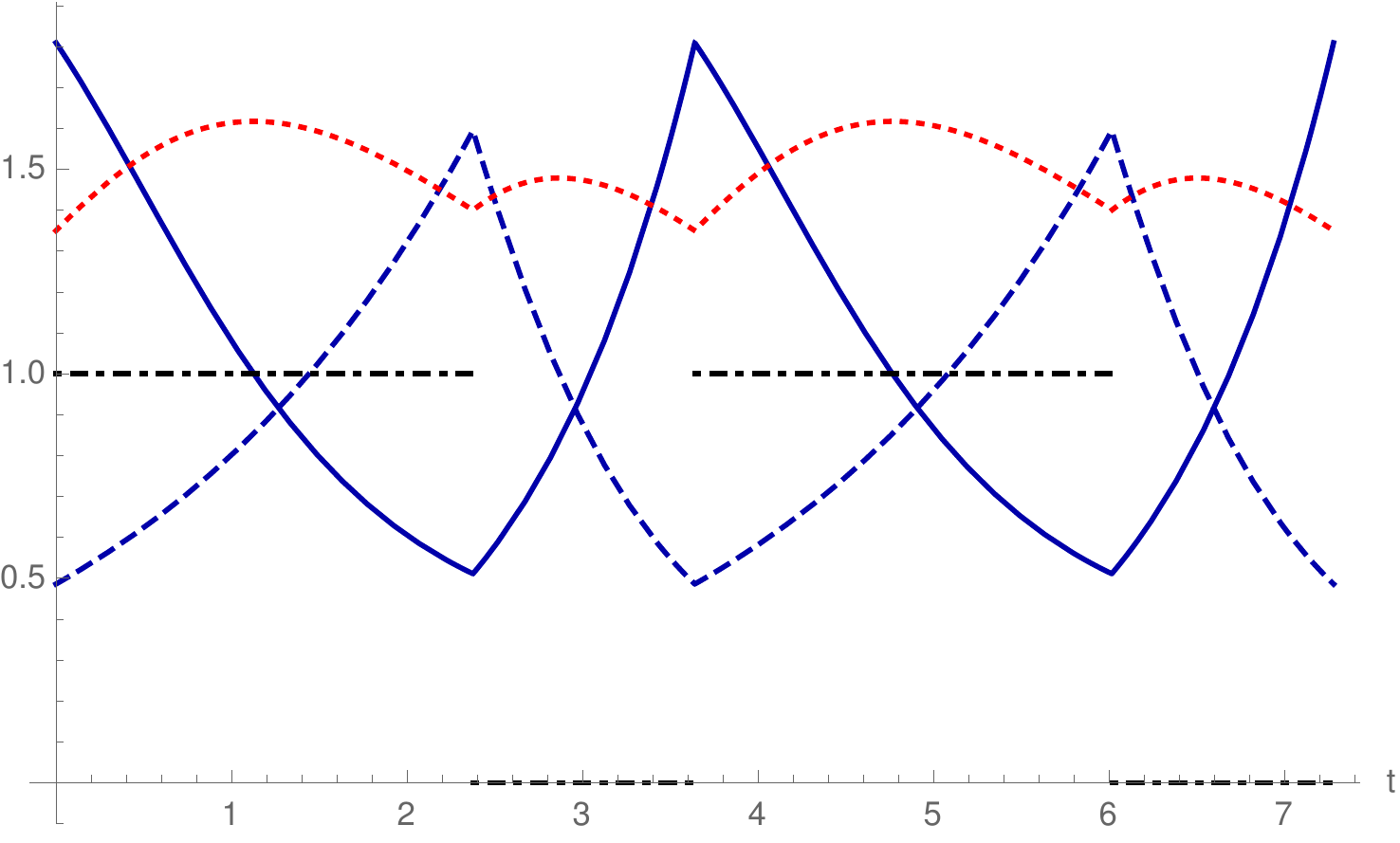}
\caption{Rescaled abundances of the preferred prey $p_1$ (solid blue curve), alternative prey $p_2$ (dashed blue curve), the predator $z$ (dotted red curve), and predator trait $q$ (dot-dashed black lines) as a function of the rescaled time $t$ (on the  horizontal axis) for a singular orbit exhibiting neither clockwise nor anticlockwise behavior (i.e., the peak in the predator density occurs almost exactly in between the two prey peaks) in the system in \eqref{rescaled_system} with $r=0.5$, $m=0.4$, and $\eps=0$. The jump points are located at $(p_1^A,p_2^A,z^A) \approx (1.81,0.49,1.35)$ and $(p_1^B,p_2^B,z^B) \approx (0.51,1.59,1.40)$.}\label{figure:appendix_fig_hybridphase_mid}
\end{figure}


\section{{Numerical continuation {of the singular periodic orbits}}}\label{section_numericalContinuation}

In this section, we use direct numerical simulations of the model system \eqref{rescaled_system} to illustrate our theoretical results from Section \ref{section_constructionOfThePeriodicOrbitSection}. The goal of this section is to highlight the role of the small parameter $\eps$. We demonstrate that we can numerically find 
\textcolor{black}{approximations to} the singular periodic orbits that we constructed in Section \ref{section_constructionOfThePeriodicOrbitSection}, which we proved to exist for `sufficiently small' $\eps$ (see Result \ref{thm:existence}). We also demonstrate that these \textcolor{black}{approximate} periodic orbits persist as $\eps$ is increased to larger values (even ones for which we no longer have a theoretical guarantee that such an orbit exists). In other words, we perform a numerical continuation in $\eps$ starting from a singular periodic orbit in which $\eps=0$.

To initiate the numerical-continuation procedure, we use the explicit analytical solution of the singular periodic orbit, which is characterized by the slow coordinates of the jump points $A$ and $B$. For a specific choice of $(r,m)$, we find the values of $(p_1^A,p_2^A,z^A)$ and $(p_1^B,p_2^B,z^B)$ that satisfy the existence conditions \eqref{AB_H0condition}, \eqref{AB_H1condition}, \eqref{AB_intp1condition}, and \eqref{AB_intp2condition}. We consider a singular periodic orbit for which a peak in the predator population lies between the peaks in the two prey populations (for a singular orbit of this kind, see Figure \ref{figure:appendix_fig_hybridphase_mid}), and we use the parameter values $(r,m) = (0.5,0.4)$. After several numerical simulations, we obtain $(p_1,p_2,z,q)\approx(1.18,0.87,1.50,0.99)$ and use these values as an initial condition, 
simulate the system \eqref{rescaled_system} with $\eps=0.025$, and find a numerical solution that is nearby the 
corresponding singular orbit (see panel (a) of Figure \ref{epsIncreasingFigure}). We carry out the continuation of this solution for increasing values of $\eps$ as follows. We simulate the rescaled system \eqref{rescaled_system} (for 50 time units when $0.025\leq\eps\leq0.5$ and for 30 time units when $0.5\leq\eps\leq1$) and use the final value of each simulation as an initial value for the next simulation in 3 sequences of 10 simulations with $\eps$ linearly spaced between $0.025$ and $0.2$ (for the simulation with $\eps=0.2$, which we show in panel (b) of Figure \ref{epsIncreasingFigure}), between $0.2$ and $0.5$ (for the simulation with $\eps=0.5$, in which we show in panel (c) of Figure \ref{epsIncreasingFigure}), and between $0.5$ and $1$ (for the simulation with $\eps=1$, which we show in panel (d) of Figure \ref{epsIncreasingFigure}).

One can clearly observe that, as the value of $\eps$ increases, the transition in $q$ between its minimum and maximum values is increasingly gradual. Moreover, {in our numerical simulations, we see that the theoretical minimal value of $q$ (i.e., $q=0$) is not attained by the solutions in panels (b), (c), and (d) of Figure \ref{epsIncreasingFigure}.
From numerical continuation, we see that choosing the small parameter $\eps = 0.025$ allows us to numerically find periodic orbits that are $\mathcal{O}(\eps)$ close to the singular limit. That is, for the parameter choice $(r,m) = (0.5,0.4)$, the value $\eps = 0.025$ can be construed as `sufficiently small' in Result \ref{thm:existence}. Additionally, we see that certain quantitative features of the singular periodic orbit persist as $\eps$ is increased. These include the antiphase oscillation between the two prey species and the occurrence of the predator peak between the prey peaks. Thus, the analytical results that we obtained through geometric singular perturbation theory, which we established for sufficiently small values of $\eps$ (see Result \ref{thm:existence}), are also meaningful for `unreasonably large' values of $\eps$. \phantomsection\label{sentence:unreasonably}

\begin{figure}[ht!]
\centering
\subfloat[$\eps=0.025$]{\includegraphics[width=0.5\textwidth]{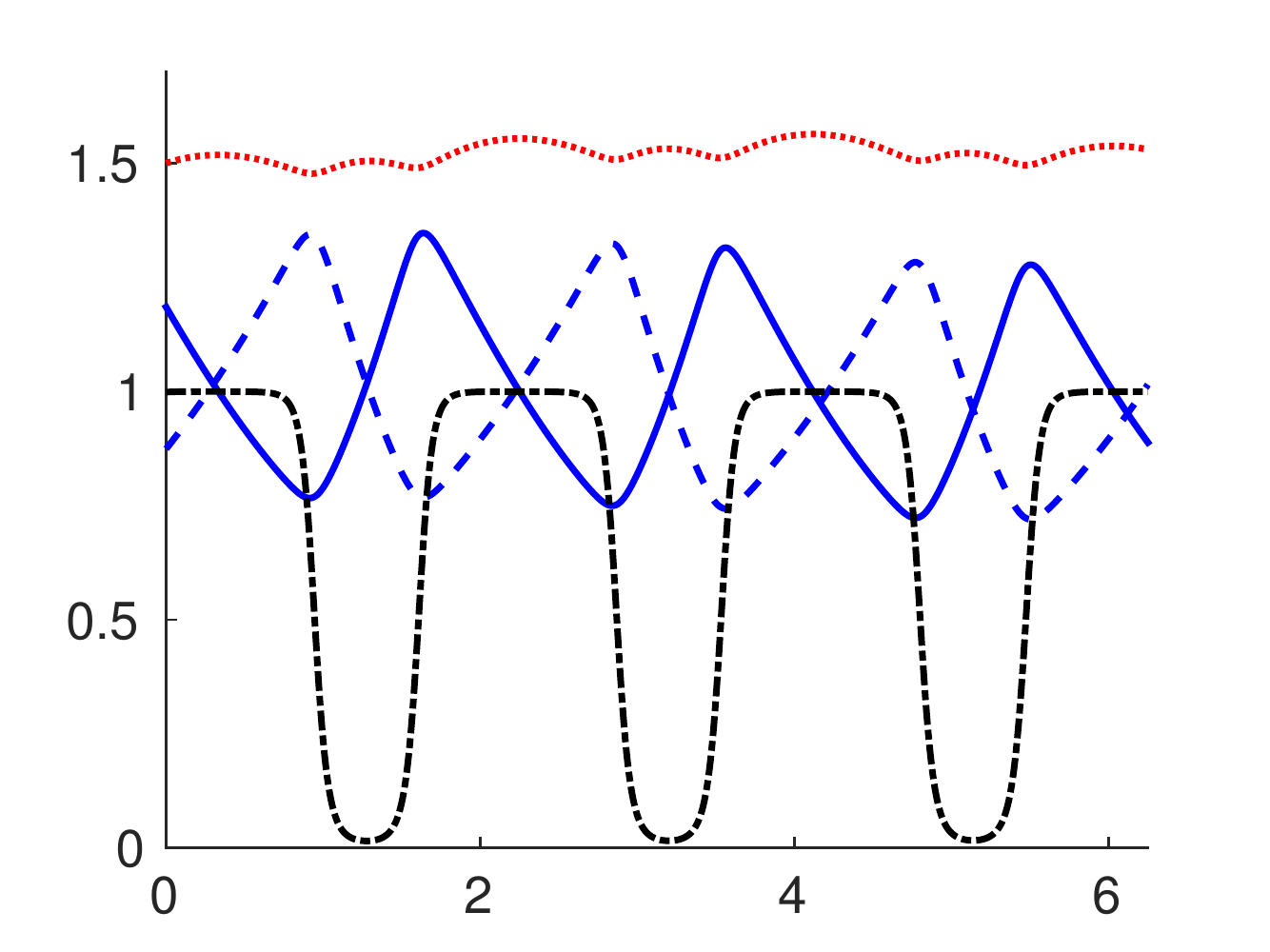}}
\subfloat[$\eps=0.2$] {\includegraphics[width=0.5\textwidth]{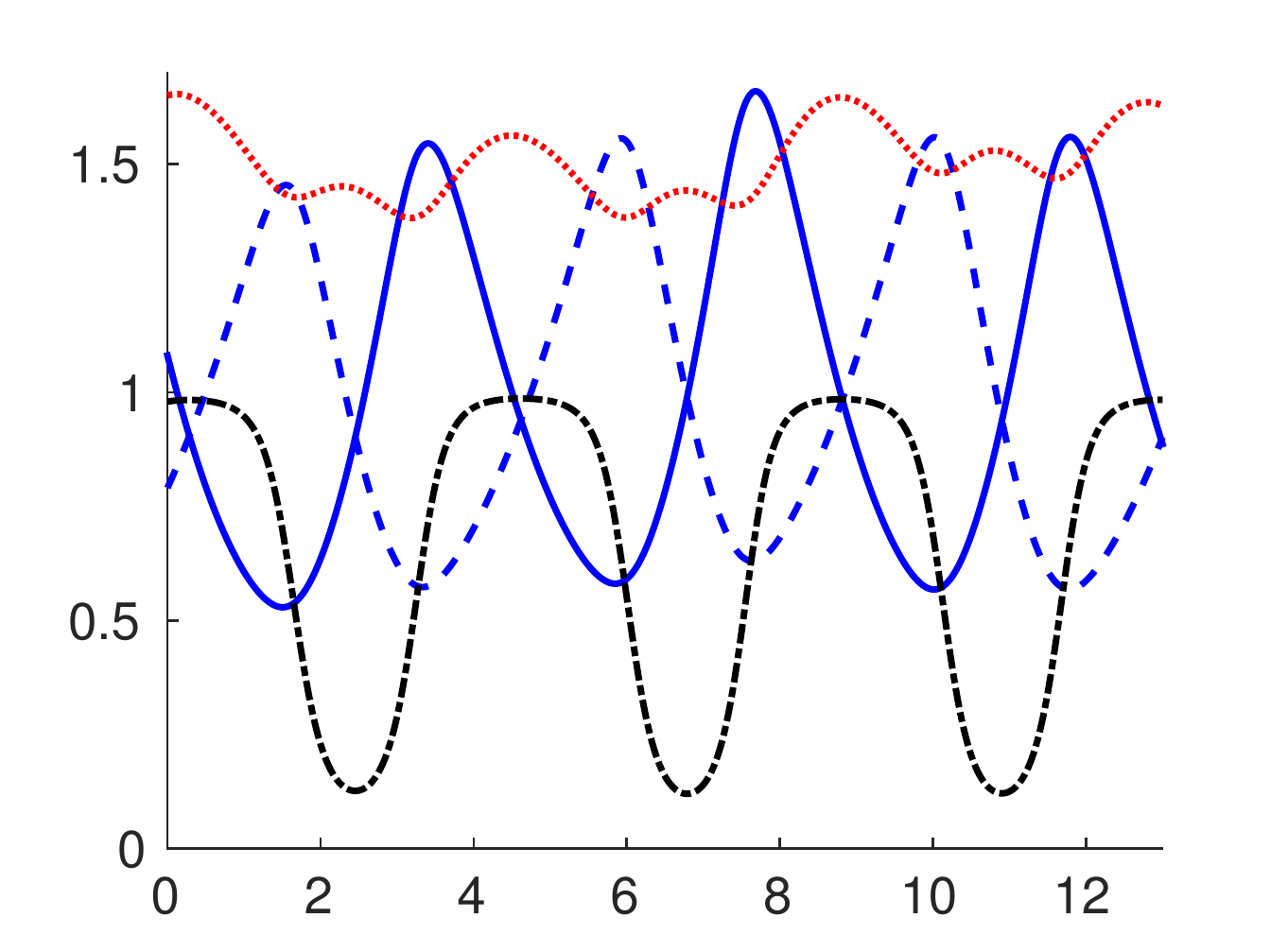}}\\
\subfloat[$\eps=0.5$] {\includegraphics[width=0.5\textwidth]{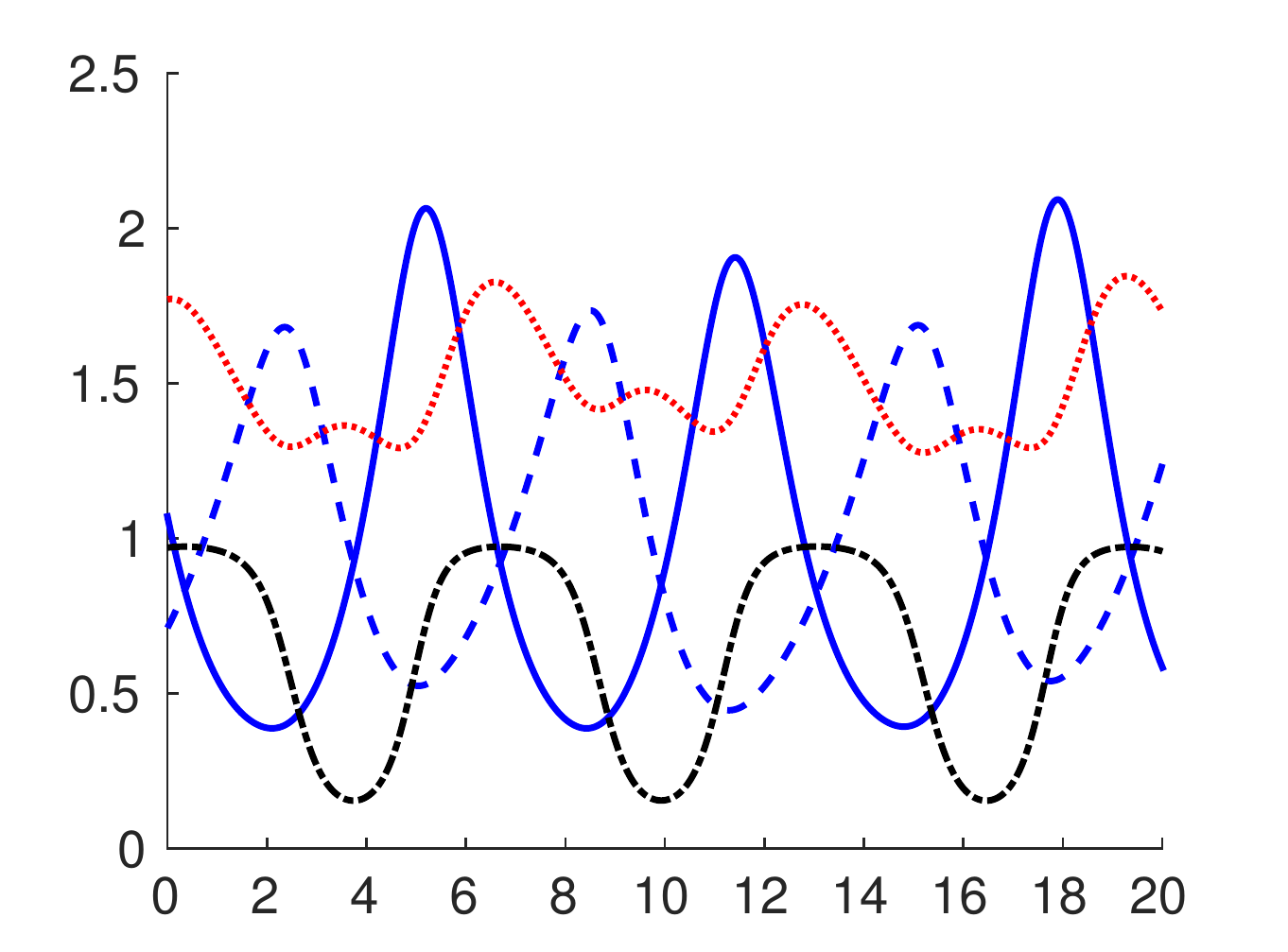}}
\subfloat[$\eps=1$]{\includegraphics[width=0.5\textwidth]{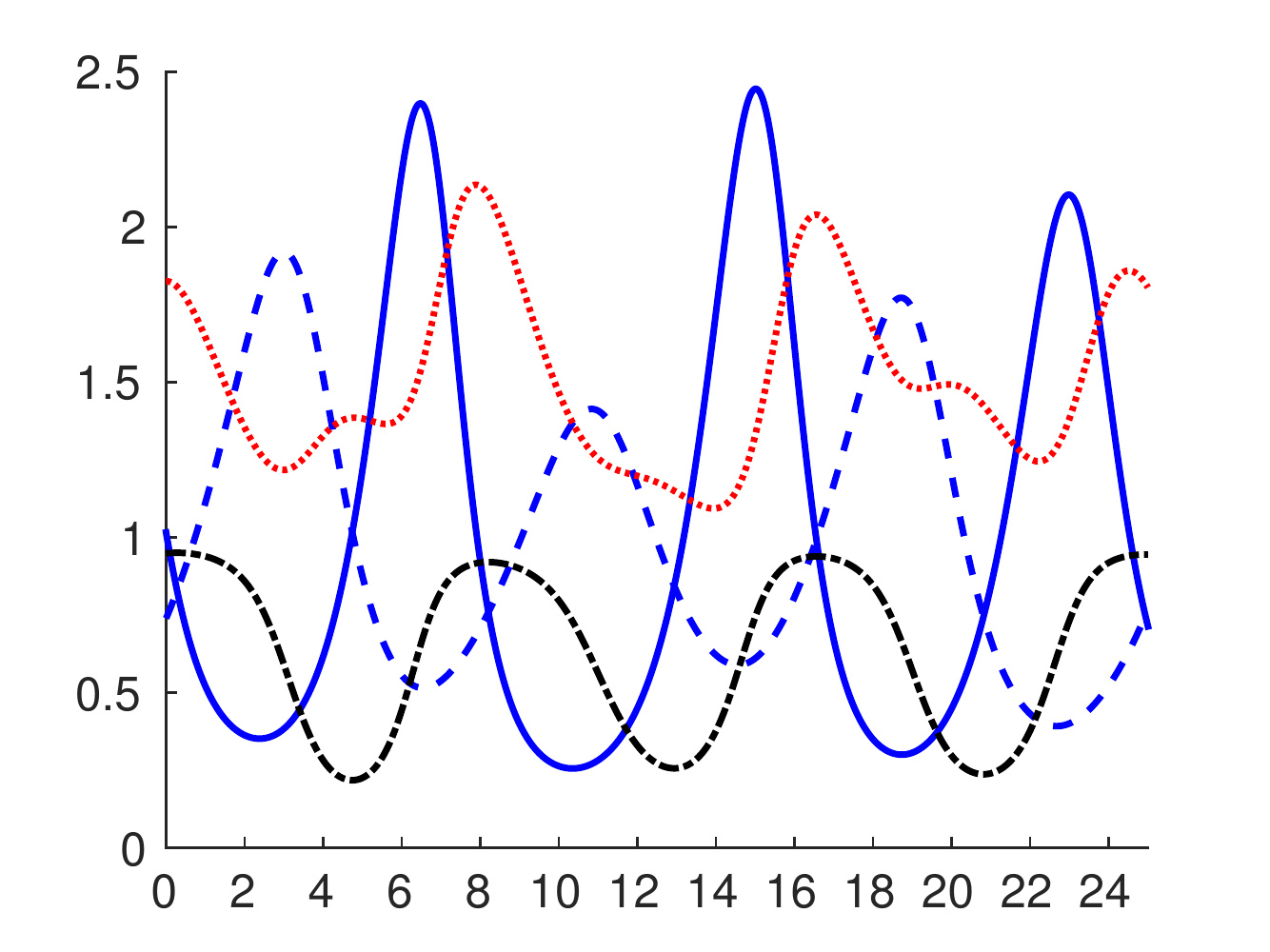}}
\caption{
Rescaled abundances of the preferred prey $p_1$ (solid blue curve), alternative prey $p_2$ (dashed blue curve), the predator $z$ (dotted red curve), and predator trait $q$ (dot-dashed black curve) as a function of the rescaled time $t$ (on the  horizontal axis) for simulations of the system in \eqref{rescaled_system} with $r=0.5$, $m=0.4$, and $\eps>0$. \textcolor{black}{We show the associated singular solution in Figure \ref{figure:appendix_fig_hybridphase_mid}.}
} 
\label{epsIncreasingFigure} 
\end{figure}


\section{Conclusions and discussion}\label{section_discussion}


\paragraph{Summary and time scales.} 

We have modeled adaptive feeding behavior of a predator that switches between two prey types. In our model (see \eqref{fastslowFullSystem}), we assumed that the predator gradually changes its diet from one prey to another depending on the prey densities and that the predator feeds only on one prey type at a time at the extremes. The change of diet is continuous, but it is fast compared to the time scale of population dynamics, so we introduced a time-scale difference between the dynamics of a predator trait (which represents the predator's desire to consume each prey) and the population dynamics. The resulting 1 fast--3 slow dynamical system exhibits periodic orbits that we can constructed analytically when the parameter $\eps$ that represents the separation of time scales is equal to $0$, and we \textcolor{black}{showed} that these orbits \textcolor{black}{can be used to find approximate periodic solutions} for small but nonzero values of $\eps$. We also demonstrated, using numerical computations, that \textcolor{black}{such approximate} periodic solutions persist for (non-small) values of $\eps$ from $\eps=0.025$ up to $\eps = 1$. 


\paragraph{Beyond piecewise-smooth formulations: Various ways to smoothen a jump.}

Part of our motivation to construct a fast--slow dynamical system for a predator adaptively switching between two prey comes from the desire to relax the assumption of a `discontinuous' predator that was made in earlier work \cite{OurPaper}. The discontinuity in this previous model, which successfully reproduces the periodicity in the ratio between the two prey groups exhibited in data collected from freshwater plankton \cite{OurPaper}, comes from the assumption that a predator chooses a diet that maximizes its growth \cite{Stephens1986}. Because it is not clear whether there exist predators that switch their feeding strategy instantaneously, our future work includes two other ways in which we `smooth out' the 1 predator--2 prey piecewise-smooth system from \cite{ourSmoothPaper}. One can regularize a piecewise-smooth dynamical system into a singular perturbation problem by `blowing up' the discontinuity boundary. This method was developed originally in \cite{Sotomayoretal1996} and later surveyed in \cite{Teixeira2012}. In the present paper, we combined these ideas with the ecological concept of fitness-gradient dynamics \cite{Lande1982,Abramsetal1993} to relate this discontinuity smoothening to a biological phenomenon. Fitness-gradient dynamics was used previously to represent trait dynamics of a predator--prey interaction in a fast--slow system for studying rapid evolution and ecological dynamics \cite{CortezEllner2010}. Note, however, that the fast--slow system in \cite{CortezEllner2010} \textcolor{black}{was} not obtained as a result of regularizing a given piecewise-smooth system. Instead, the rate of change of a predator (or prey) trait was assumed to be governed by fitness-gradient dynamics \cite{Lande1982,Abramsetal1993} and \textcolor{black}{was assumed} to evolve at a faster time scale than that of the predator--prey interaction \cite{CortezEllner2010}\textcolor{black}{; we have used both of these assumptions in the present paper.} 


\paragraph{Obtaining analytical results using time-scale separation.}

An important feature of our model \eqref{rescaled_system} is the presence of the parameter $\eps$, which introduces a time-scale separation between the (fast) model component $q$ and the (slow) model components $(p_1,p_2,z$). The inclusion of this time-scale separation enables us to not only include a biological mechanism (i.e., natural selection) directly into the model}, but also to {prove the existence of \textcolor{black}{approximate} periodic orbits for sufficiently small values of this time-scale separation parameter $\eps$. Using numerical calculations, we demonstrated that these `sufficiently small' values of $\eps$ are within numerical reach, as we are able to find values for $\eps$ for which direct numerical simulation of the model \eqref{rescaled_system} yields periodic orbits that are $\mathcal{O}(\eps)$ close to their singular counterparts. By using numerical continuation, we have also shown that these orbits persist for increasing values of $\eps$. Therefore, the method that we outlined in this paper can be used not only to study periodic solutions to \eqref{rescaled_system} in the presence of a time-scale separation, but also (using numerical continuation of the singular periodic solution) to investigate the existence and behavior of periodic solutions to \eqref{rescaled_system} when the time scale of the rate of change of the trait $q$ is \emph{comparable} to that of the predator--prey interaction. The periodic solutions that we constructed in this paper are \emph{far-from-equilibrium} solutions, so the model variables do not stay close to the system's coexistence equilibrium. As we mentioned in Section \ref{section_lin_centrecentre}, this equilibrium is of center--center type. In general, the use of local analysis around this equilibrium to study periodic solutions is extremely complicated and intricate, and the nature of this type of analysis excludes the study of far-from-equilibrium solutions (and, in particular, it excludes ones in which the trait variable $q$ switches between $0$ and $1$.}

\textcolor{black}{As we stated in Remark \ref{rmk:noSotoTrevino}, the presence of only one fast component prohibits the use of `standard' existence results from the literature on geometric singular perturbation theory. In particular, the number of slow directions is related directly to the lack of transversality of the intersection of the stable and unstable manifolds of the invariant manifolds $\mathcal{M}_0$ and $\mathcal{M}_1$. The existence problem would have to be unfolded to higher orders in $\eps$ to obtain a subset of singular periodic solutions that persist for all time as fully periodic solutions. Numerical simulations of the }\textcolor{black}{system \eqref{rescaled_system}} \textcolor{black}{indeed show that not every approximate periodic solution remains bounded for long times, and several numerical periodic solutions exhibit a slowly modulated amplitude.}
\textcolor{black}{These phenomena, and the problem of `true' persistence of periodic orbits, are interesting subjects for future research.}

It is worth noting that the method that we employed in this paper to construct (singular) periodic orbits is not confined to orbits with two slow segments (one on each slow hyperplane $\mathcal{M}_{0,1}$, as in Figure \ref{figure:per_fastred}). 
Using the same methods, our analysis can be extended to study periodic orbits with two slow segments on each slow hyperplane by concatenating them using four fast transitions. This would lead to an extension of the family of possible periodic orbits. \textcolor{black}{This larger class of periodic orbits, which exhibit a wider range of qualitative features, can also be fit to experimental data. For more discussion on comparison with experimental data, see the last paragraph of this discussion section.}


\paragraph{Rapid evolution versus phenotypic plasticity.}

The two mechanisms of adaptivity---i.e., phenotypic plasticity and rapid evolution---cause rapid adaptation and affect population dynamics \cite{Shimadaetal2010,Yamamichietal2011}. Although it is not clear precisely how these different mechanisms affect population dynamics, it has been suggested that models that account for phenotypic plasticity exhibit a stable equilibrium more often than models that account for rapid evolution \cite{Yamamichietal2011}. It has also been suggested that this situation can arise from a faster response time of plastic genotypes than that of nonplastic genotypes to fluctuating environmental conditions \cite{Yamamichietal2011}. Indeed, our model \eqref{fastslowFullSystem}, which describes the population dynamics of a predator and its two prey in the presence of rapid evolutionary change in a predator trait, does not contain stable steady states.
In contrast, the model in \cite{OurPaper}, which considers an adaptive change of diet in response to prey abundance, exhibits convergence to a steady state for a large parameter range. 


\paragraph{Cryptic and out-of-phase cycles.}

Empirical evidence suggests that rapid evolution is a possible mechanistic explanation for cyclic dynamics that differ from those in traditional predator--prey systems. For an evolving prey, such dynamics include (1) large-amplitude cycles in a predator population while a prey population remains nearly constant \cite{yoshida2007cryptic,BohannanLenski2000}, (2) predator and prey oscillating almost exactly out of phase \cite{Becks2010,Yoshidaetal2003}, and (3) oscillations in which a peak in a prey population follows that in a predator population \cite{Becks2010}. Oscillations of types (1) and (2) also arise in predator--prey models with rapid predator evolution \cite{CortezEllner2010}. It has also been demonstrated that rapid predator evolution as a response to an evolving prey can generate cyclic dynamics both in experiments \cite{HiltunenBecksNatComm2014} and in models \cite{CortezTheorEcol2015} of coevolution. Our model \eqref{fastslowFullSystem} represents an evolving predator that feeds on two different types of prey and exhibits different types of periodic orbits---including ones in which the predator and prey populations oscillate out of phase, total prey density remains approximately constant, and a peak in the prey population follows that in the predator populations. 


\paragraph{Future work and comparison with experiments and field observations.}

To identify orbits that exist in an ecologically reasonable parameter range, our ongoing work includes comparing our model simulations with data collected from freshwater plankton in the field \cite{TirokGaedke2006,TirokGaedke2007a}. Our principal model assumptions (i.e., large population size, short generation times, and well-mixed environment) hold for these data. However, similar models can be formulated for any other organisms that satisfy these assumptions, including the microorganisms used in the laboratory experiments in \cite{HiltunenBecksNatComm2014}. The insight into rapid evolution gained from studying a tractable plankton system can be used as an example for understanding rapid evolution of larger organisms and their abilities to adapt to changing environmental conditions (such as climate change or species introductions). Moreover, one of the major applications of the understanding of coevolution in microorganisms is resistance to antimicrobial drugs. By fitting parameters of prey growth rates and predator mortality to data, we expect to be able to distinguish parameter regimes to determine which members of periodic-orbit families best describe the data, and one can thereby gain insights into a system of one evolving predator that feeds on two different types of prey. In particular, we expect such comparisons between models and data to help determine ecological trade-offs and their possible influence on rapid predator evolution. Empirical evidence from a study of coevolving predator and prey suggests that a predator pays a low fitness cost (or no cost at all) for counteracting anti-predatory prey evolution \cite{HiltunenBecksNatComm2014}. However, in addition to the unknown mechanism of the predator response, the trade-off(s) that constrain rapid predator evolution remain unknown if one only looks at data without doing any modeling.


\section*{Acknowledgements} 

We thank Stephen Ellner, Christian Kuehn, and David J. B. Lloyd for helpful discussions. We thank Ursula Gaedke for sending us data. We also thank the anonymous referees and the editor for their constructive feedback and comments that helped to improve this paper. The Lake Constance data were obtained within the Collaborative Programme SFB 248 funded by the German Science Foundation. SHP was supported by Osk. Huttunen Foundation, Engineering and Physical Sciences Research Council through the Oxford Life Sciences Interface Doctoral Training Centre, and People Programme (Marie Curie Actions) of the European Union's Seventh Framework Programme (FP7/2007-2013) under REA grant agreement \#609405 (COFUNDPostdocDTU). FV was supported by \href{http://www.nwo.nl/en}{NWO} through a Rubicon grant.


\bibliography{fastSlow1pred2prey_siadsManuscript6}{}
\bibliographystyle{siam}


\appendix
\section{Geometric singular perturbation theory}\label{appendix_gspt}

\textcolor{black}{In the present paper, we gained insight into how the evolution of traits that occur on a comparable time scale to that of ecological interactions arises in population dynamics by studying the limit in which trait evolution occurs on a much faster time scale than that of the predator--prey interactions. Consequently, our goal in choosing a geometric approach was to create solution trajectories for a parameter $\eps>0$ by concatenating segments of curves that are determined by either the fast reduced dynamics or the slow reduced dynamics when $\eps=0$. In this appendix, we give a brief introduction to this kind of procedure. See \cite{Jones1995,KuehnMTS2015,Hek2010} for further details.}

Following the notation in \cite{KuehnMTS2015}, a fast--slow dynamical system with $m$ fast variables and $n$ slow variables (and time as the only independent variable) is expressed as
\begin{align}
	\eps \frac{\text{d}  x}{\text{d} t}&=\eps \dot{x}=f(x,y,\eps)\,,\nonumber\\
	\frac{\text{d} y}{\text{d} t}&=\dot{y}=g(x,y,\eps)\,,
\label{system_in_slowTime}
\end{align}
where $f:\,\mathbb{R}^m \times \mathbb{R}^n \times \mathbb{R} \to \mathbb{R}^m$, $g:\,\mathbb{R}^m \times \mathbb{R}^n \times \mathbb{R} \to \mathbb{R}^n$, and $\eps$ (with $0< \eps \ll1$) is the ratio of the two time scales. We rescale the slow time $t$ by $\eps$ and obtain an equivalent system that evolves on the fast time scale $\tau=t/\eps$. We thus write
\begin{align}
	\frac{\text{d} x}{\text{d} \tau}&=x'=f(x,y,\eps)\,,\nonumber\\
	\frac{\text{d} y}{\text{d} \tau}&=y'=\eps g(x,y,\eps)\,.
\label{system_in_fastTime}
\end{align} 
We can take the (singular) limit $\eps\rightarrow0$ in \eqref{system_in_slowTime}, which describes the dynamics evolving on the slow time scale $t$, to obtain the reduced slow vector field
\begin{align}
	0&=f(x,y,0)\,,\nonumber\\
	\dot{y}&=g(x,y,0)\,.
\label{slow_reduced_problem}
\end{align}
We call \eqref{slow_reduced_problem} the {\emph{slow reduced problem}}. Similarly, we can take the limit $\eps\rightarrow0$ in \eqref{system_in_fastTime}, which describes the dynamics evolving on the fast time scale $\tau$, to obtain the reduced fast vector field
\begin{align}
	x'&=f(x,y,0)\,,\nonumber\\
	y'&=0\,,
\label{fast_reduced_problem}
\end{align}
which is called the {\emph{fast reduced problem}}.
The fast and slow reduced problems are connected through the {\emph{critical manifold}} $C_0=\{(x,y) \in \mathbb{R}^m \times \mathbb{R}^n: f(x,y,0)=0 \}$, a sufficiently smooth submanifold of $\mathbb{R}^m \times \mathbb{R}^n$. The critical manifold $C_0$ determines the equilibrium points of the fast reduced problem \eqref{fast_reduced_problem}, and the differential-algebraic slow reduced problem \eqref{slow_reduced_problem} determines a (slow) dynamical system on $C_0$. 

By using geometric singular perturbation theory, one studies the critical manifold and the slow \eqref{slow_reduced_problem} and fast \eqref{fast_reduced_problem} reduced problems to obtain information about the behavior of the full system \eqref{system_in_slowTime}. This approach builds on the work of Fenichel \cite{Fenichel1979}, which guarantees that, under some general conditions, several geometric objects (e.g., the critical manifold $C_0$) defined in the reduced slow and fast problems persist for sufficiently small $\eps > 0$ as similar geometric objects in the full system \eqref{system_in_slowTime}. Hence, for example, (singular) orbits that are constructed by concatenating orbit pieces from the slow reduced problem \eqref{slow_reduced_problem} and the fast reduced problem \eqref{fast_reduced_problem} persist for sufficiently small $\eps>0$ in the sense that such a singular orbit is a good approximation of a `true' orbit of the full system \eqref{system_in_slowTime}. 

For an introduction to geometric singular perturbation theory and its concepts, see \cite{Hek2010}. For a comprehensive overview of singular perturbation theory, see \cite{KuehnMTS2015}.


\section{Finding solution families}\label{appendix:findsols}

To find singular periodic orbits to \eqref{rescaled_system}, we seek values for $(p_1^A,p_2^A,z^A)$ and $(p_1^B,p_2^B,z^B)$ that satisfy equations \eqref{AB_H0condition}, \eqref{AB_H1condition}, \eqref{AB_intp1condition}, and \eqref{AB_intp2condition}. We highlight aspects of the procedure by analyzing these equations in detail, starting with \eqref{AB_intp1condition}.

We obtain the integral in \eqref{AB_intp1condition} by integrating the slow dynamics on $\mathcal{M}_1$ [see equation \eqref{T1_direct}]. As we mentioned in Section {\ref{section_analysis_of_theFastSlowSystem}}, the slow coordinates $(p_1,z)$ form a Lotka--Volterra system (see \eqref{slowFlow_onqIsOne}). Because all orbits in the $(p_1,z)$-system are closed, the function $z(p_1)$ needed for the integrand of \eqref{AB_intp1condition} cannot be determined uniquely. Indeed, a full (closed) Lotka--Volterra orbit in the $(p_1,z)$ system consists of two branches of $z(p_1)$. As we illustrated in Figure \ref{figure:LVonM1}, there is both an upper branch and a lower branch. Using \eqref{p1_func_z} to solve $z(p_1)$, we {use the Lambert \emph{W}-function $W_i(x)$ \cite{dlmf}} to explicitly write the two branches as
\begin{align}\label{eq:z0-1_p1}
	 {z_+\left(p_1;p_1^A,z^A\right)} &
	 = -W_{-1}\left(-z^A e^{-z^A} \left(\frac{p_1^A e^{-p_1^A}}{p_1 e^{-p_1}}\right)^m \right)\,,
	 \\
	 {z_-\left(p_1;p_1^A,z^A\right)} &
	 = -W_{0}\left(-z^A e^{-z^A} \left(\frac{p_1^A e^{-p_1^A}}{p_1 e^{-p_1}}\right)^m \right)\,,
\end{align}
where ${z_+}$ indicates the upper branch and ${z_-}$ indicates the lower branch of the associated Lotka--Volterra orbit. The branches connect at the left and right extrema, which are given respectively by $(p_1^\text{min},1)$ and $(p_1^\text{max},1)$; again see Figure \ref{figure:LVonM1}).

\begin{figure}[t]
 \centering
\includegraphics[width=0.6\textwidth]{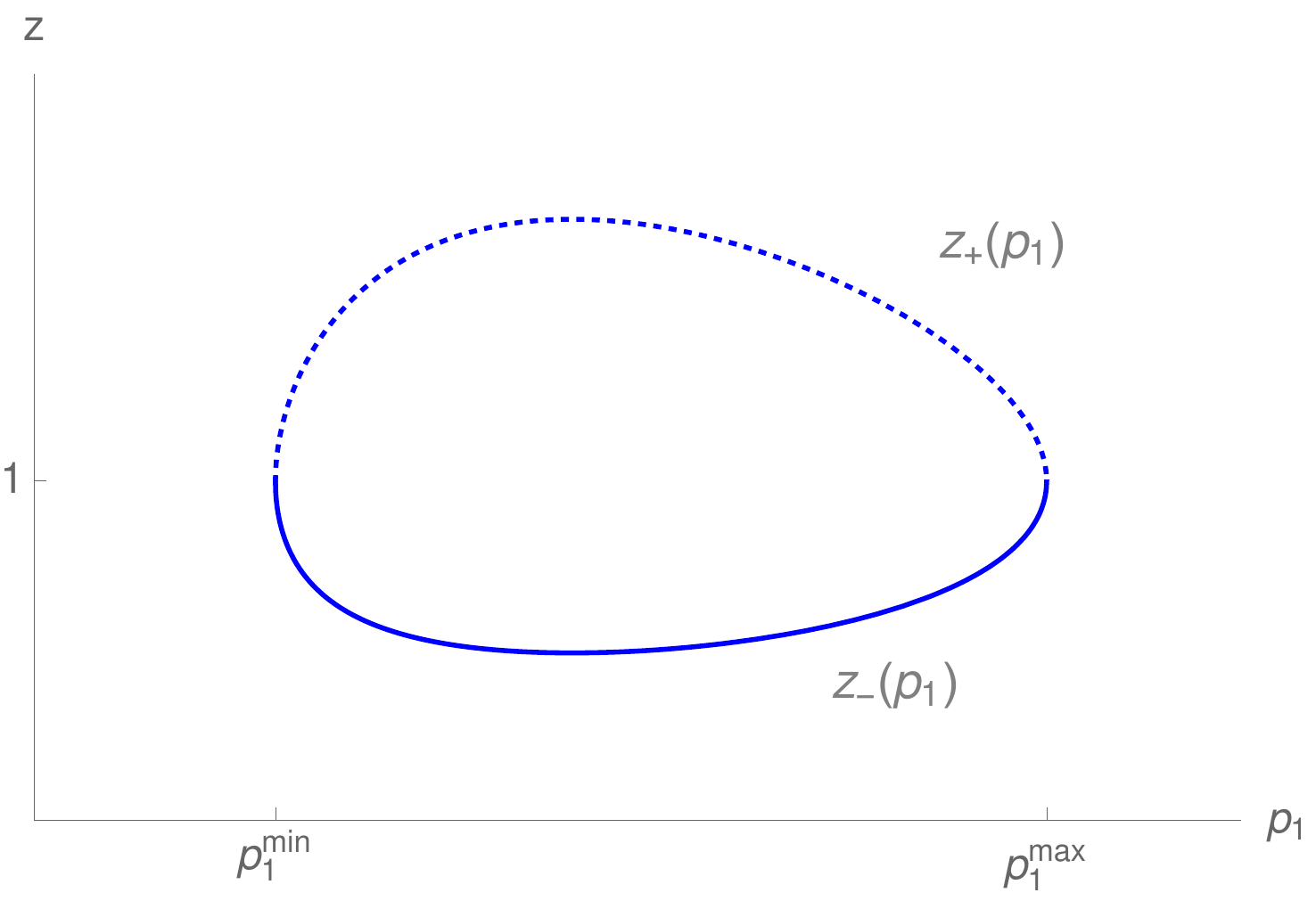}
\caption{{A Lotka--Volterra orbit in the $(p_1,z)$ phase plane. The orbit is closed and consists of two branches, $z_+$ (dotted blue curve) and $z_-$ (solid blue curve); see \eqref{eq:z0-1_p1}. The branches connect at the left and right extrema, which are located at $(p_1^\text{min},1)$ and $(p_1^\text{max},1)$, respectively; see \eqref{eq:p1minmax}.}}\label{figure:LVonM1}
\end{figure}

The final expression of the integrand \eqref{AB_intp1condition} depends on the path followed by the orbit on $\mathcal{M}_1$. In Figure \ref{figure:M1examples}, we show two examples of such paths. In the first example, the initial point $(p_1^A,z^A)$ is in the upper left quadrant (i.e., $p_1^A<1$ and $z^A>1$). The final point $(p_1^B,z^B)$ is in the lower left quadrant (i.e., $p_1^B<1$ and $z^B<1$). For this path, the slow travel time $T_1$ \eqref{T1_direct} is given by
\begin{equation}\label{eq:M1example1_T1}
 	\int_{p_1^A}^{p_1^\text{min}} \frac{1}{1-{z_+}(p_1)}\,\frac{\text{d} p_1}{p_1} + \int_{p_1^\text{min}}^{p_1^B} \frac{1}{1-{z_-}(p_1)}\,\frac{\text{d} p_1}{p_1}\,.
\end{equation}
In the second example, the initial point lies in the lower right quadrant (i.e., $p_1^A>1$ and $z^A<1$). The slow travel time $T_1$ is then
\begin{equation}\label{eq:M1example2_T1}
	 \int_{p_1^A}^{p_1^\text{max}} \frac{1}{1-{z_-}(p_1)}\,\frac{\text{d} p_1}{p_1} + \int_{p_1^\text{max}}^{p_1^\text{min}} \frac{1}{1-{z_+}(p_1)}\,\frac{\text{d} p_1}{p_1} + \int_{p_1^\text{min}}^{p_1^B} \frac{1}{1-{z_-}(p_1)}\,\frac{\text{d} p_1}{p_1}\,.
\end{equation}

\begin{figure}[t]
 \centering
\includegraphics[width=0.49\textwidth]{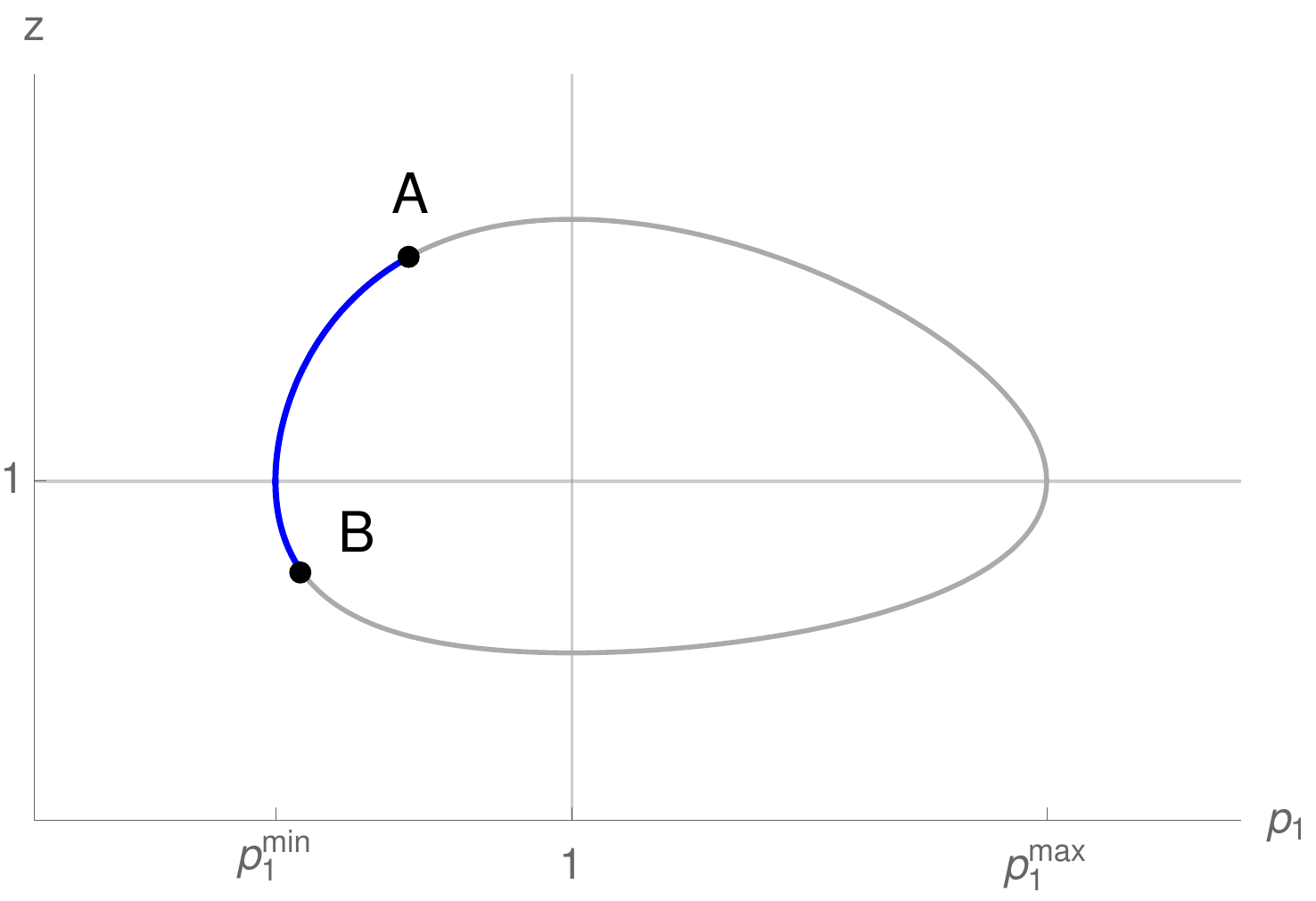}
\includegraphics[width=0.49\textwidth]{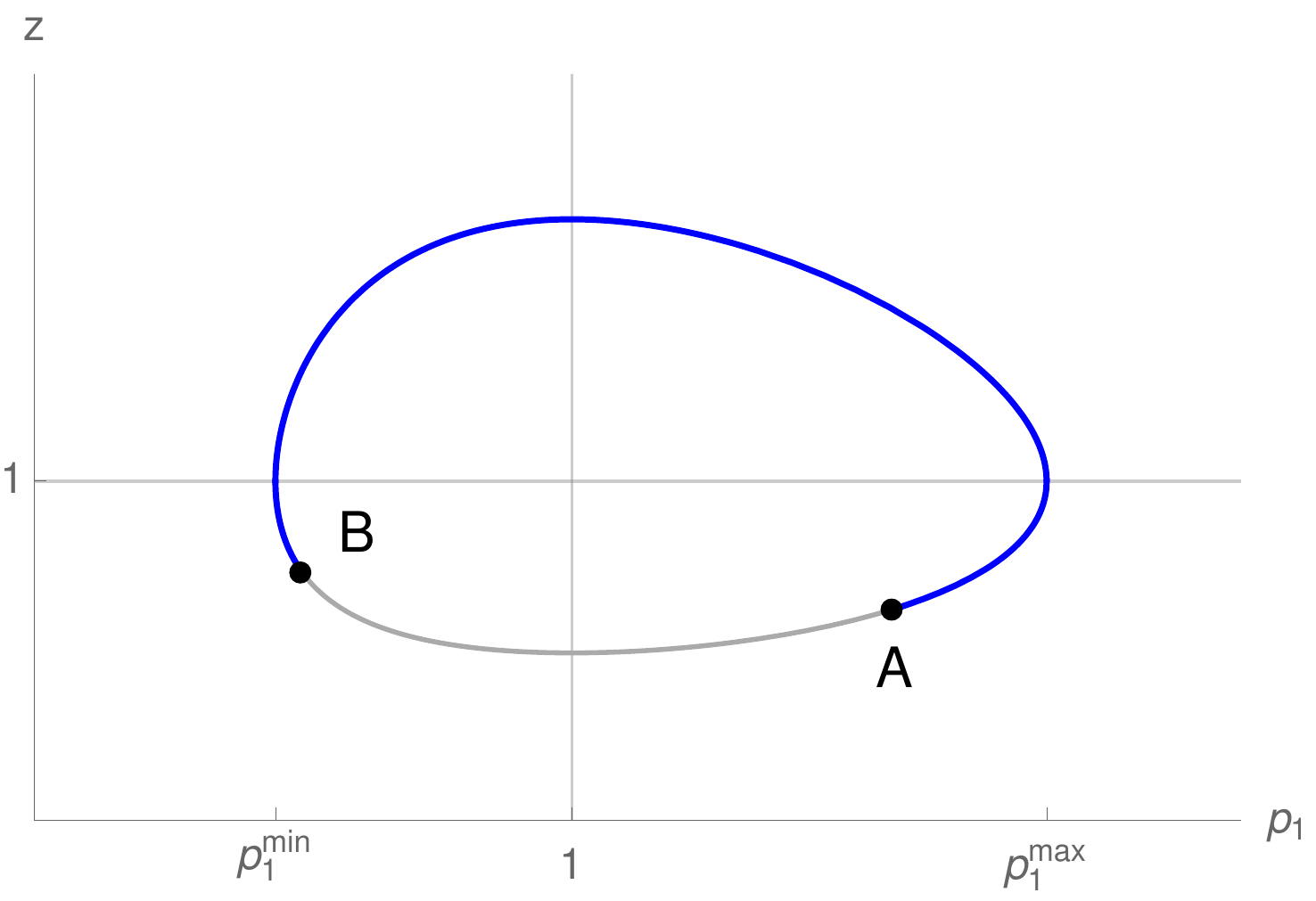}
\caption{Example paths of the slow dynamics on $\mathcal{M}_1$. The travel time $T_1$ for the left path is given by \eqref{eq:M1example1_T1}, and the travel time for the right path is given by \eqref{eq:M1example2_T1}.}\label{figure:M1examples}
\end{figure}

Taking into account all possible combinations of initial and final points, a methodical analysis yields the following explicit expression for condition \eqref{AB_intp1condition}:
\begin{equation}\label{AB_intp1condition_explicit}
	\frac{1}{r} \log \left(\frac{p_2^B}{p_2^A}\right) = \int_{p_1^A}^{p_1^B}\frac{1}{1-{z_+}(p_1)} \frac{\text{d}p_1}{p_1} + \left\{\begin{aligned} \int_{p_1^\text{min}}^{p_1^B}\frac{1}{1-{z_-}(p_1)}-\frac{1}{1-{z_+}(p_1)} \frac{\text{d}p_1}{p_1}\,, & &\text{if }z^B<1\,, \\ \int_{p_1^A}^{p_1^\text{max}}\frac{1}{1-{z_-}(p_1)}-\frac{1}{1-{z_+}(p_1)} \frac{\text{d}p_1}{p_1}\,, & &\text{if }z^A<1\,, \\ 0\,, & &\text{otherwise\,,} \end{aligned}\right.
\end{equation}
where ${z_\pm}$ are defined in \eqref{eq:z0-1_p1}. One can use the {Lambert \emph{W}-function} to give explicit expressions for the extremal values $p_1^\text{min,max}${, yielding}
\begin{equation}\label{eq:p1minmax}
\begin{aligned}
	 p_1^\text{min}(p_1^A,z^A) &= -W_0\left(-p_1^A e^{-p_1^A} \left(z^A e^{1-z^A}\right)^\frac{1}{m}\right)\,,\\
	 p_1^\text{max}(p_1^A,z^A) &= -W_{-1}\left(-p_1^A e^{-p_1^A} \left(z^A e^{1-z^A}\right)^\frac{1}{m}\right)\,.
\end{aligned}
\end{equation}

A detailed analysis of \eqref{AB_intp2condition} on $\mathcal{M}_0$ {is analogous to that of} \eqref{AB_intp1condition}. {Similar to \eqref{eq:z0-1_p1}, we introduce}
\begin{align}\label{z0-1_p2}
	{\zeta_{+}\left(p_2;p_2^A,z^A\right) }&{= -r\,W_{-1}\left(-\frac{z^A}{r} e^{-\frac{z^A}{r}} \left(\frac{p_2^A e^{-p_2^A}}{p_2 e^{-p_2}}\right)^\frac{m}{r} \right)\,,}\\
	 {\zeta_{-}\left(p_2;p_2^A,z^A\right) }&{= -r\,W_{0}\left(-\frac{z^A}{r} e^{-\frac{z^A}{r}} \left(\frac{p_2^A e^{-p_2^A}}{p_2 e^{-p_2}}\right)^\frac{m}{r} \right)}
\end{align}
to explicitly express condition \eqref{AB_intp2condition} as
\begin{equation}\label{AB_intp2condition_explicit}
	 \log \left(\frac{p_1^A}{p_1^B}\right) = \int_{p_2^B}^{p_2^A}\frac{1}{r-{\zeta_+}(p_2)} \frac{\text{d}p_2}{p_2} + \left\{\begin{aligned} \int_{p_2^\text{min}}^{p_2^A}\frac{1}{r-{\zeta_-}(p_2)}-\frac{1}{r-{\zeta_+}(p_2)} \frac{\text{d}p_2}{p_2}\,, & &\text{if }z^A<r\,, \\ \int_{p_2^B}^{p_2^\text{max}}\frac{1}{r-{\zeta_-}(p_2)}-\frac{1}{r-{\zeta_+}(p_2)} \frac{\text{d}p_2}{p_2}\,, & &\text{if }z^B<r\,, \\ 0\,, & &\text{otherwise\,,} \end{aligned}\right.
\end{equation}
where the extremal values $p_2^\text{min,max}$ are given by
\begin{equation}
\begin{aligned}
 	p_2^\text{min}(p_2^A,z^A) &= -W_0\left(-p_2^A e^{-p_2^A}\left(\frac{z^A}{r} e^{1-\frac{z^A}{r}}\right)^{\frac{r}{m}}\right)\,,\\
 	p_2^\text{max}(p_2^A,z^A) &= -W_{-1}\left(-p_2^A e^{-p_2^A}\left(\frac{z^A}{r} e^{1-\frac{z^A}{r}}\right)^{\frac{r}{m}}\right)\,.
\end{aligned}
\end{equation}

Finally, we use the conserved-quantity conditions \eqref{AB_H0condition} and \eqref{AB_H1condition} to eliminate $p_2^B$ and $p_1^B$ from \eqref{AB_intp1condition_explicit} and \eqref{AB_intp2condition_explicit}. Again using the {Lambert \emph{W}-function}, we obtain
\begin{align}
 	p_2^B &= -W_{0,-1}\left(-p_2^A e^{-p_2^A} \left(\frac{z^A}{z^B} e^{\frac{z^B-z^A}{r}}\right)^{\frac{r}{m}}\right)\,,	\label{p2B_explicit}\\
 	p_1^B &= -W_{0,-1}\left(-p_1^A e^{-p_1^A} \left(\frac{z^A e^{-z^A}}{z^B e^{-z^B}}\right)^{\frac{1}{m}}\right)\,.\label{p1B_explicit}
\end{align}
Substituting \eqref{p2B_explicit} and \eqref{p1B_explicit} into \eqref{AB_H0condition} and \eqref{AB_H1condition} yields a pair of rather lengthy conditions that must be satisfied by $(p_1^A,p_2^A,z^A)$ (i.e., the three slow coordinates of $A_{0,1}$) and $z^B$. In Figure \ref{figure:exconds_sample}, we show numerical computations of these solutions for {the parameter values $(r,m) = (0.5,0.4)$}.

\begin{figure}[t]
 \centering
\includegraphics[width=0.49\textwidth]{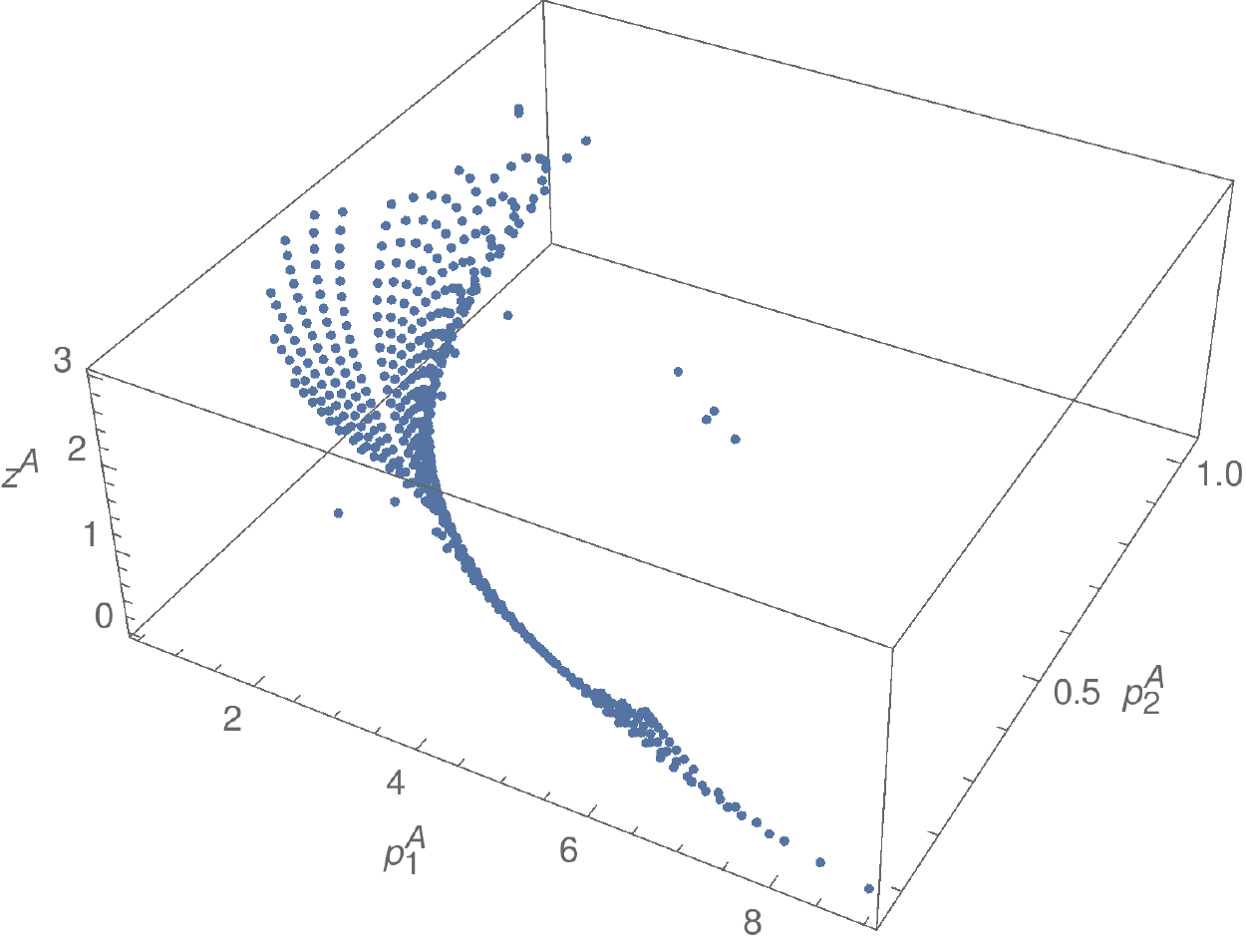}
\includegraphics[width=0.49\textwidth]{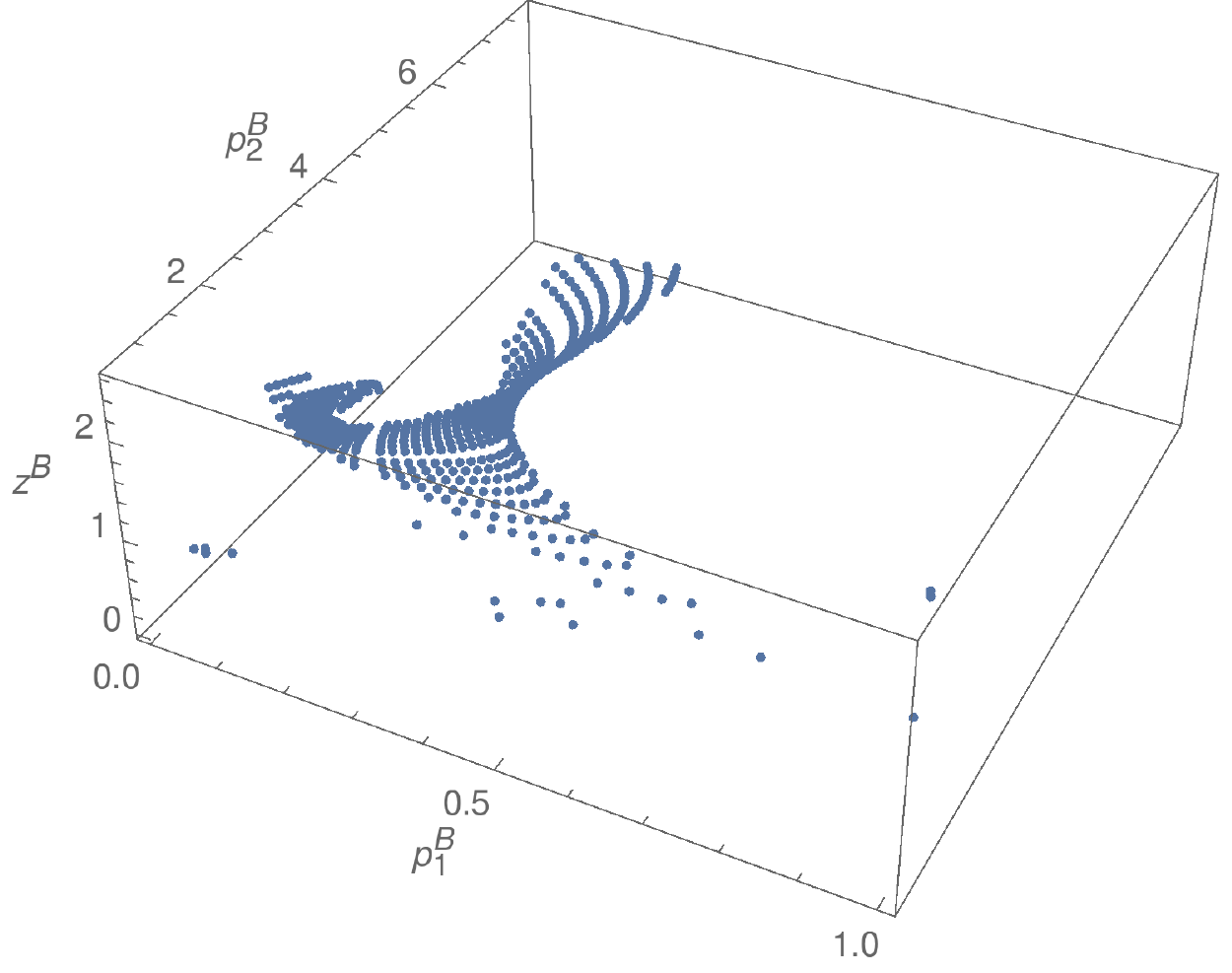}
\includegraphics[width=\textwidth]{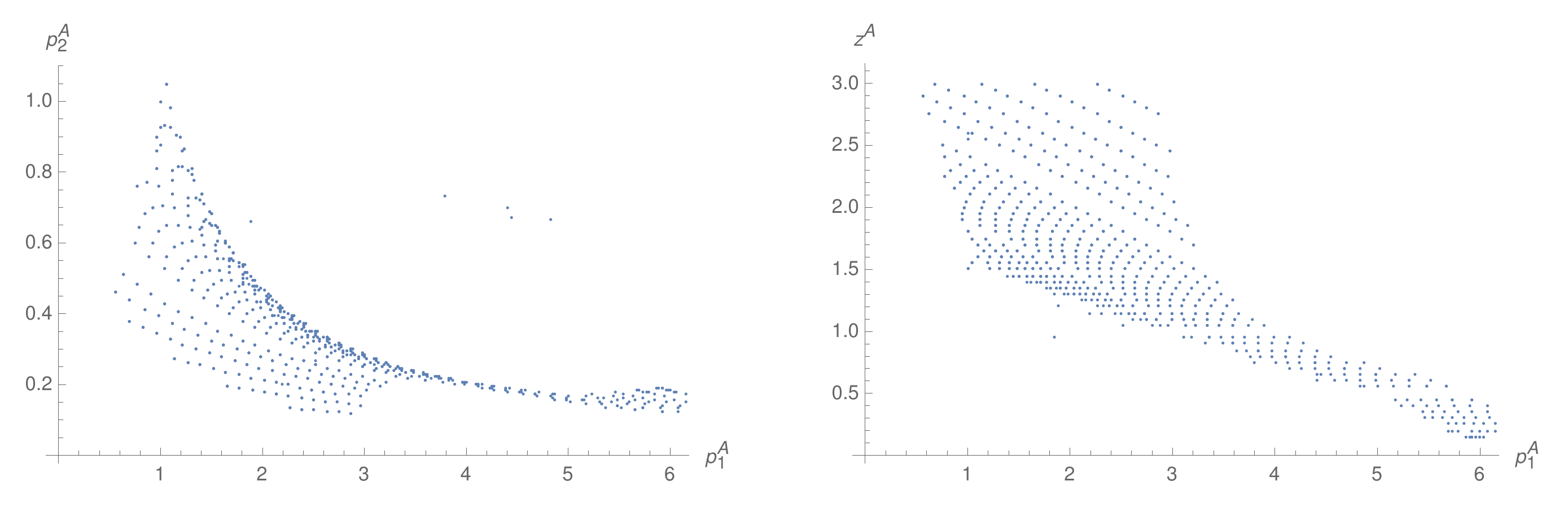}
\includegraphics[width=\textwidth]{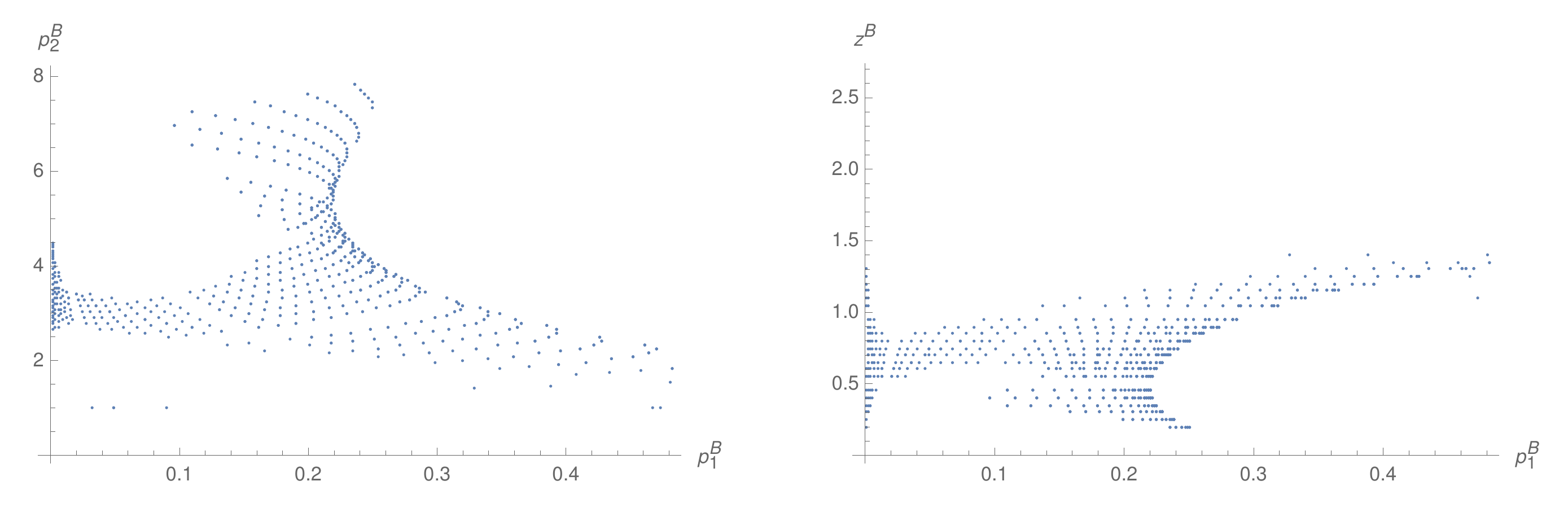}
\caption{\textcolor{black}{The coordinates of $(p_1^A,p_2^A,z^A)$ and $(p_1^B,p_2^B,z^B)$ for singular periodic orbits with prey--prey synchronization for $(r,m) = (0.5,0.4)$. Each dot represents a numerical solution of the existence conditions \eqref{AB_H0condition}, \eqref{AB_H1condition}, \eqref{AB_intp1condition}, and \eqref{AB_intp2condition} using the explicit formulation in Appendix \ref{appendix:findsols}. Each such numerical solution is thus given as a six-tuple $((p_1^A,p_2^A,z^A),(p_1^B,p_2^B,z^B))$. In the panels of this figure, we show the projections of these six-tuples onto different coordinate planes.}}\label{figure:exconds_sample}
\end{figure}

\end{document}